  \documentstyle[epsfig, 11pt]{article}

\newcommand{\qed}{{$\Box$}}

\newcommand{\ba}{\begin{array}}
\newcommand{\ea}{\end{array}}

\newtheorem{lemma}{Lemma}
\newtheorem{hypo}{Assumption}
\newcommand{\bhypo}{\begin{hypo}}
\newcommand{\ehypo}{\end{hypo}}
\newtheorem{defi}{Definition}
\newcommand{\ble}{\begin{lemma}}
\newcommand{\ele}{\end{lemma}}
\newcommand{\bde}{\begin{defi}}
\newcommand{\ede}{\end{defi}}

\newtheorem{prop}{Proposition}
\newcommand{\epr}{\end{prop}}
\newcommand{\bpr}{\begin{prop}}
\newtheorem{teo}{Theorem}
\newcommand{\bth}{\begin{teo}}
\newcommand{\eth}{\end{teo}}
\newtheorem{rema}{Remark}
\newcommand{\bre}{\begin{rema}}
\newcommand{\ere}{\end{rema}}
\newcommand{\ee}{\end{equation}}
\newcommand{\be}{\begin{equation}}

\begin{document}

\centerline{\Large\bf  Pole Distribution of PVI Transcendents
 } 
\vskip 0.3 cm
\centerline{\Large\bf  close to a Critical Point (summer 2011)} 
\vskip 0.3 cm
\centerline{\Large Davide Guzzetti}
\vskip 0.2 cm
\centerline{SISSA, Intenational School of Advanced Studies, 34136 Trieste, Italy}
\vskip 0.2 cm 
\centerline{and KIAS, Korea Institute of Advanced Study, 130-722 Seoul, South Korea}
\vskip 0.2 cm
\centerline{Tel +82 -2-958-3861}
\vskip 0.2 cm
\centerline{davide$_{-}$guzzetti@yahoo.com}

\begin{abstract}   

\vskip 0.2 cm
 The distribution of the poles of  Painlev\'e VI transcendents associated to  semi-simple Frobenius manifolds is determined  close to a critical point.  It is shown that the poles accumulate at the critical point, asymptotically along two rays. 
As an  example, the Frobenius manifold given by the quantum cohomology of $CP^2$ is considered.  The general PVI is also considered. 
\end{abstract}

\section{Introduction}

 Consider the Sixth Painlev\'e equation associated to a three dimensional Frobenius manifold, hereafter denoted $PVI_{\mu}$, with parameters  $\alpha={(2\mu-1)^2\over 2}$, $\beta=\gamma=0$, $\delta={1\over 2}$ (in standard notation). 
$$
{d^2y \over dx^2}={1\over 2}\left[ 
{1\over y}+{1\over y-1}+{1\over y-x}
\right]
           \left({dy\over dx}\right)^2
-\left[
{1\over x}+{1\over x-1}+{1\over y-x}
\right]{dy \over dx}
$$
$$
+{1\over 2}
{y(y-1)(y-x)\over x^2 (x-1)^2}
\left[
(2\mu-1)^2+{x(x-1)\over (y-x)^2}
\right]
,~~~~~\mu \in {\bf C}.~~~~~~~~~~~\hbox{\rm ($PVI_\mu$)}
$$

The algebraic solutions of  $PVI_\mu$ were studied in \cite{DM}, its elliptic representation in \cite{D4}. The importance of $PVI_\mu$ in the theory of semi--simple Frobenius manifolds was extabilished in \cite{Dub1} and \cite{Dub}, and practically applied in \cite{guzz1} to the construction of some relevant manifolds.  

\vskip 0.2 cm
In this paper, we study the distribution of the movable poles of Painlev\'e transcendents,  close to the critical point  $x=0$. We do this for any $\mu\in {\bf C}$. Due to the symmetries of $PVI_\mu$, the results, obtained close to $x=0$, can be translated to the poles close to $x=1$ and $x=\infty$. We shown that the poles accumulate at the critical point, asymptotically along two rays. The results can be extended to the general Painlev\'e VI equation, as sketched in section \ref{generalPVI}. 

\vskip 0.2 cm
 The distribution of the poles close to a critical point for $PVI_\mu$ has been anticipated in \cite{D4}, where it was  conjectured that the poles of a transcendent, considered as the meromorphic extension of a branch on the universal covering of the critical point, should accumulate at the critical point {\it along  spirals}.  The same conjecture is motivated in \cite{D3} for the general Painlev\'e VI equation. 
 In \cite{Br}, the pole distribution for  PVI with parameters $\alpha=\beta=\gamma=1/8$ and $\delta=3/8$ (Hitchin's equation \cite{hitchin}) is determined  on the whole universal covering of
 ${\bf C}\backslash\{0,1,\infty\}$. A formula for an infinite series of poles is given in terms of Theta-functions. The poles are distributed along  lines which  are spirals at a small scale  around the critical points, and more complicated lines on the whole universal covering. A birational Okamoto's transformation  transforms the above PVI equation into $PVI_{\mu}$ with $\mu=1/2$. The latter, studied by Picard \cite{Picard}, will be considered in section \ref{exempicard}.

\vskip 0.2 cm 
 In this paper, we determine the poles   of Painlev\'e transcendents in a neighborhood of $x=0$, with
 bounded $\arg x$, namely  $|\arg x|<\vartheta$, for some $\vartheta>0$. Thus, $x$ may tend to zero along a
 {\it radial path}, while spiral paths
 are not allowed.  With this limitation, "most" solutions of PVI have no poles in a sufficiently small
 neighborhood of the critical point, except for a class of solutions, which is the object of this paper. 
 Actually, the critical behaviors at one of the critical points $x=0,1$ or $\infty$ of a branch $y(x)$ can
 be divided into a few classes, classified in \cite{guz2010} and \cite{guztable} (an equivalent
 classification is given in \cite{Bruno6} \cite{Bruno7}). 
 Among them,  one class of three-real parameters
 solutions   {\it may admit poles in a neighborhood of the critical point, with bounded $\arg x$}.  The class is
 given in Proposition 1 of \cite{guz2010} for the general  Painlev\'e VI equation. In the particlular case
 of $PVI_\mu$, the proposition states the following:

\vskip 0.2 cm 

Let 
$\nu \in {\bf R}$, $\nu\neq 0$, and  $d\in {\bf C}$ be given, such that $2\mu-1\neq \pm 2i\nu$. Let also $\vartheta>0$ be given.  Equation $PVI_\mu$ has a solution $y(x)$  admitting the following expansion when $x\to 0$ and $|\arg x|<\vartheta$: 
\be
\label{pro2y}
{1\over y(x)} = y_1(x)+xy_2(x)+x^2y_3(x)+...~=\sum_{n=1}^\infty x^{n-1} y_n(x)
\ee
where: 
\be
\label{coefA}
 y_n(x)=\sum_{m=-n}^n A_{nm}(\nu,\mu)~ e^{2imd}x^{2im\nu}
\ee
For any $\vartheta$, there exists a sufficiently small $\epsilon$ such that the series of $1/y(x)$ converges in the domain $0<|x|<\epsilon$, $|\arg x|<\vartheta$, and defines  an holomorphic function of $x$ and $x^{2i\nu}$ . 
The $A_{nm}(\nu,\mu)$'s are   rational functions of $\nu$, $\mu$, and satisfy the property $\bar{A}_{n,m}(\nu,\mu)= A_{n,-m}(\nu,\mu)=A_{nm}(-\nu,\mu)$ 
 (the bar denotes the complex conjugate). 
 Their explicit form is recursively   computed by the procedure of \cite{guz2010}. For example, the lower order coefficients are: 
$$
A_{11}={(2\mu-1-2i\nu)^2\over 16\nu^2},~~~~~A_{10}= -{(2\mu-1)^2-4\nu^2\over 8\nu^2}.
$$
and the first order approximation is: 
\be
\label{firstordermu}
y_1(x)=
          {(2\mu-1+2i\nu)^2e^{-2id}x^{-2i\nu}\over 16 \nu^2}-{(2\mu-1)^2-4\nu^2\over 8\nu^2}
+
{(2\mu-1-2i\nu)^2e^{2id}x^{2i\nu}\over 16\nu^2}
\ee
$$
\equiv 
1-{4\nu^2+(2\mu-1)^2\over 4 \nu^2} \sin^2\left( \nu \ln x +d+{i\over 2}\ln\left[{2\mu-1+2i\nu\over 2\mu-1-2i\nu}\right]\right)
$$
The second order coefficients are:
$$ 
A_{22}= -{(-2 \mu + 2 i \nu + 1)^4 \over  2^9 \nu^4 },~~~
A_{21}= {( (2\mu-1)^2   + 8 i \nu^3 )(-2 \mu + 2 i \nu + 1)^2 \over 2^7 \nu^4},
$$
$$
A_{20}= -{ ( (2\mu-1)^2  + 4 \nu^2) (3(2\mu^2  - 1)^2 - 4 \nu^4)\over 2^8 \nu^4}.
$$
 The other coefficients up to order $y_4$ are in the Appendix of the 
preprint version arXiv:1104.5066. 
Due to the structure of $y_n(x)$, which is invariant for $\nu\mapsto -\nu$ and $d\mapsto d+k\pi$, $k$ integer,  we are allowed to assume that 
\be
\label{Nishigamo}
 0\leq \Re d \leq \pi,~~~~~\nu>0.
\ee
We remark that (\ref{pro2y}) is derived in \cite{guz2010} by a symmetry transformation applied to a solution of a PVI equation with $\alpha=\gamma=1-2\delta=0$ and $\beta\neq 0$,  defined for $\arg x$ bounded, namely $|\arg x|<\vartheta$, for some $\vartheta>0$. The latter solution can be locally constructed both by the method of \cite{Jimbo}  or  the method of local analysis of \cite{Sh} (see also \cite{IKSY}) and \cite{D3}. 
It follows from these methods that  $\vartheta$ is chosen finite but  arbitrarily, observing that, if $\vartheta$ is increased, the radius of convergence of  (\ref{pro2y}) in general decreases.

% % % % % % % % % % % % % 

\subsection{Poles close to the Critical Point and main Results of the Paper}

Clearly, the $n$-th order  $y_n(x)$ is an oscillatory bounded function in a neighborhood of $x=0$. The leading term of $1/y(x)$ is  then $y_1(x)$, while $x^{n-1}y_n(x)=O(x^{n-1})$. 
The solution $y(x)$ with expansion (\ref{pro2y})  may have poles in a neighborhood of $x=0$, which are the zeros of $\sum_{n=1}^\infty x^{n-1} y_n(x)$. The zeros of $y_1(x)$ do not coincide in general with the zeros of $\sum_{n=1}^\infty x^{n-1} y_n(x)$. Thus, one cannot  write: 
$$
  y(x) = {1\over y_1(x)}\left[ 1-x ~{y_2(x)\over y_1(x)}+O(x^2)\right],~~~x\to 0
$$
The above asymptotic expansion is  true only when  $x\to 0$ in a sector not containing the zeros of both $y_1(x)$ and $\sum_{n=1}^\infty x^{n-1} y_n(x)$ (in this case, the leading term $1/y_1(x)$ of the expansion is computed  also in \cite{Bruno6} and  \cite{Bruno7}).  We expect that the poles of $y(x)$ are  close to the zeros of $y_1(x)$ as $x\to 0$, because $y(x)^{-1}\sim y_1(x)$. It is to be remarked that any considerations about the poles of $y(x)$ must be done for $|x|$ smaller than the radius of convergence of the series (\ref{pro2y}).

%% THEOREM %%%%%%%%%%%%%%%%%%%%%%%%%%%%%%%%%%%%%%%%%%%%%
\vskip 0.2 cm 
  The general result of the paper is summarized  in the following theorem. 

\bth
\label{theo1}
 Let 
$\nu>0$ and  $d\in {\bf C}$, $0\leq \Re d \leq \pi$,  be given such that $2\mu-1\neq \pm 2i\nu$, and let $y(x)$ be (\ref{pro2y}).  
  Then, $y_1(x)$ has the two sequences of zeros

$$
   x_k(1)= \exp\left\{-{d\over \nu}-{k \pi \over \nu}\right\}= \exp\left\{-i{\Im d\over \nu}\right\}
\exp\left\{-{\Re d\over \nu} - {k \pi \over \nu}\right\},~~~~~~~~~~~k\in{\bf Z},
$$
and
$$
   x_k(2)=x_k(1)~\exp\left\{-{i\over \nu} \ln\left({2\mu-1+2i\nu\over 2\mu-1-2i\nu}\right)\right\}=~~~~~~~~~~~~~~~~~~~~~~~~~~~~~~~~~~~~
$$
$$~~~=
 x_k(1)~\exp\left\{{1\over \nu} \arg\left({2\mu-1+2i\nu\over 2\mu-1-2i\nu}\right)\right\}~
\exp\left\{-{i\over \nu}\ln\left|{2\mu-1+2i\nu\over 2\mu-1-2i\nu}\right|\right\}.
$$ 
where  $-\pi<\arg\left({2\mu-1+2i\nu\over 2\mu-1-2i\nu}\right)\leq \pi$, being other choices  absorbed into  ${k \pi \over \nu}$.
Let $k_0\in{\bf N}$  be sufficiently big in order for $x_k(j)$, $j=1,2$ to fall in the domain of convergence of the series (\ref{pro2y}). There exists $K$ sufficientley big such that for every  $k\geq \max\{K,k_0\}$, and every $j=1,2$, $y(x)$ has a  pole $\xi_k(j)$ lying in a neighborhood of $x_k(j)$,  with   the following asymptotic representation 
\be
\label{POLESintro}
   \xi_k(j)=x_k(j)-{1\over 2} x_k(j)^2 +\sum_{N=3}^\infty\Delta_N(j) x_k(j)^N, ~~~k\to +\infty,~~~x_k(j)\to 0.
\ee
The coefficients $\Delta_N(j)\in{\bf C}$ are certain numbers independent of $k$  that can be computed form the coefficients $A_{nm}$ of (\ref{coefA}). The first terms  are:

$$
\Delta_3(1) = {16\mu^4-8\mu^2+176\nu^4+352\nu^2+177\over 1024(\nu^2+1)^2}
$$
$$
\Delta_3(2) = {16\mu^4-64\mu^3+88\mu^2-48\mu+176\nu^4+352\nu^2+185\over1024(\nu^2+1)^2}
$$
$$
\Delta_4(1) = -{16\mu^4-8\mu^2+49+48\nu^4+96\nu^2   \over 1024(\nu^2+1)^2}
$$
$$
\Delta_4(2) = -{16\mu^4-64\mu^3+88\mu^2-48\mu+57+48\nu^4+96\nu^2 \over 1024(\nu^2+1)^2}
$$
\eth

\qed

\begin{figure}
\epsfxsize=9cm
\epsfysize=9cm
\centerline{\epsffile{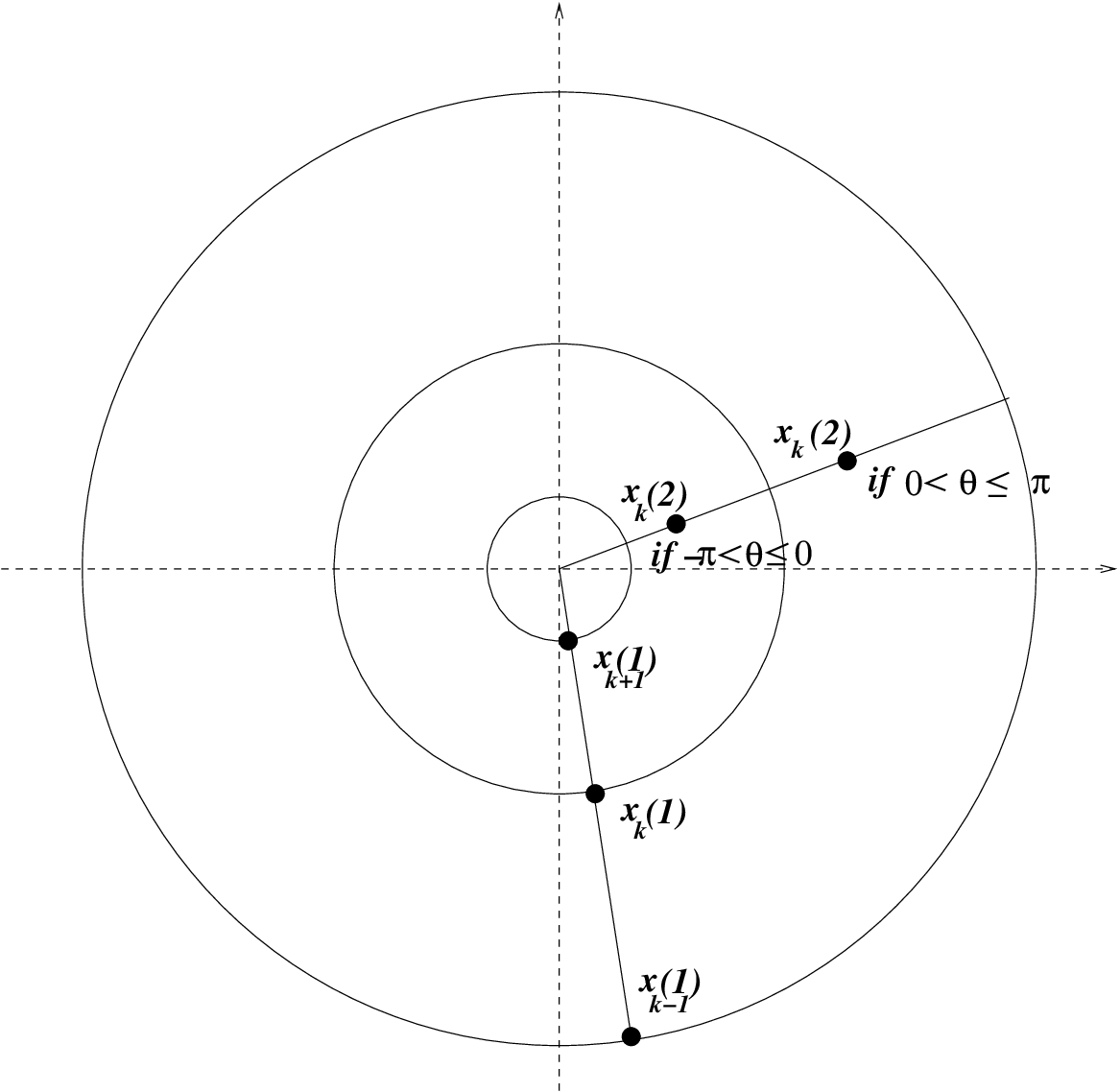}}
\caption{Position of the zeros of $y_1(x)$.  $\theta= \arg\left({2\mu-1+2i\nu\over 2\mu-1-2i\nu}\right)$. Observe that $|x_{k+1}(1)|<|x_k(2)|\leq |x_k(1)|$ when $-\pi < \arg\left({2\mu-1+2i\nu\over 2\mu-1-2i\nu}\right)\leq 0 $, and  $|x_k(1)|<|x_k(2)|\leq |x_{k-1}(1)|$ when 
 $0< \arg\left({2\mu-1+2i\nu\over 2\mu-1-2i\nu}\right)\leq \pi $.}
\label{orderzero}
\end{figure}

%%% END THEOREM %%%%%%%%%%%%%%%%%%%%%%%%%
\vskip 0.3 cm 
In the following, when the $x_k(j)$, j=1,2, are considered, it will be understood that $k\geq k_0$. 
When $k\to +\infty$, the zeros $x_k(1)$  accumulate at $x=0$   along the ray of angle $-{\Im d\over \nu}$, while the zeros $x_k(2)$  accumulate at $x=0$  along the ray of angle $-{1\over \nu}\left(\Im d + \ln\left|{2\mu-1+2i\nu\over 2\mu-1-2i\nu}\right|\right)$. A typical case is that of figure \ref{orderzero}. For real $\mu$ the two rays coincide.  If  $\mu=1/2$, then $\Delta_N(1)=\Delta_N(2)$ and the poles of the two sequences overlap (double poles).  The series (\ref{POLESintro}) is at least asymptotic, but we can prove its convergence when $\mu={1\over2}$ (section \ref{exempicard}). 
Observe that in order for the rays of the zeros  to fall into the domain where the series (\ref{pro2y}) is defined, $\vartheta$ must be sufficiently big (see also the remark at the end of section \ref{PVIFROBENIUS}).   .

\begin{figure}
\epsfxsize=9cm

\epsfysize=9cm

\centerline{\epsffile{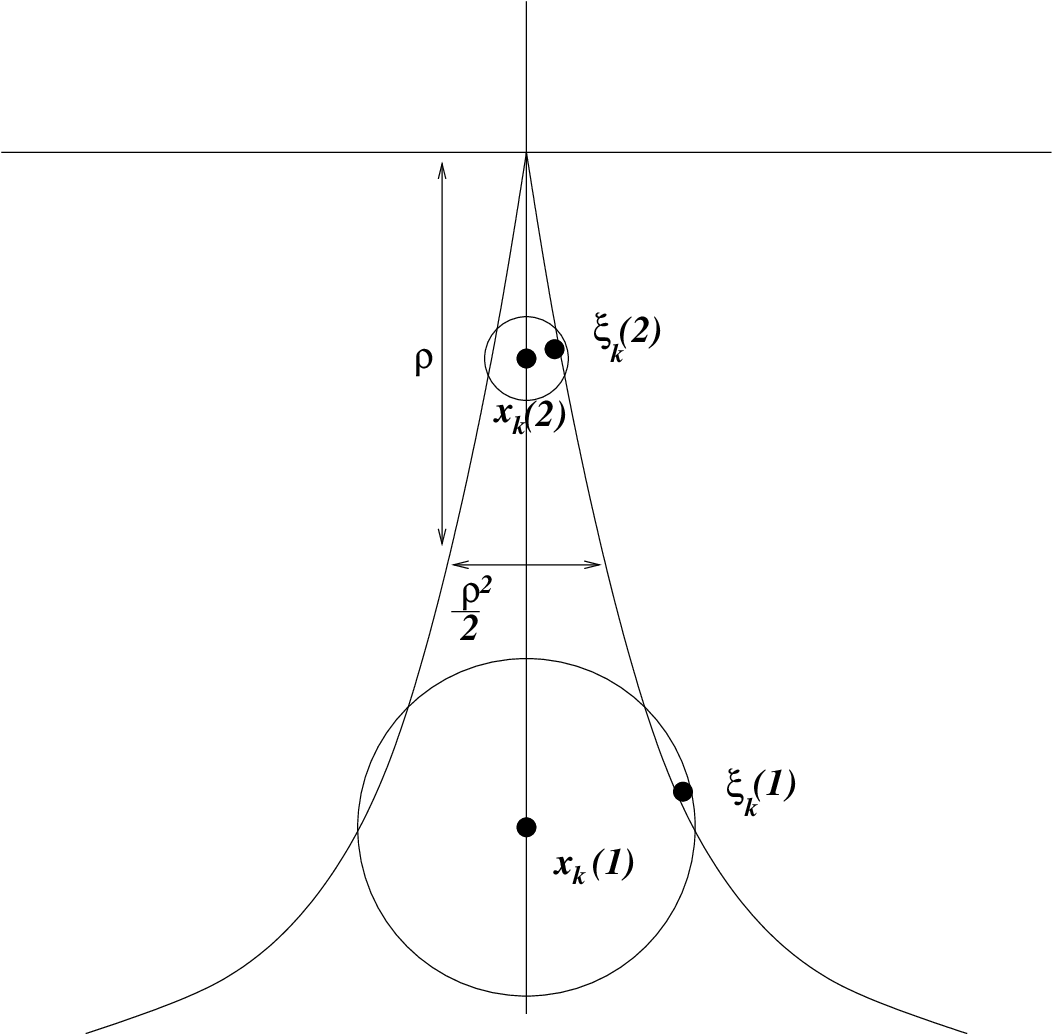}}
\caption{The disks around the zeros $x_k(j)$ of $y_1(x)$, where the poles $\xi_k(j)$  of $y(x)$ possibly lie. The figure represents the case $\mu=-1$, $\nu=2\ln {\bf G}/\pi$, ${\bf G}=(1+\sqrt{5})/2$. In this case,  the two rays coincide with the negative imaginary axis.} 
\label{poles}
\end{figure}

\vskip 0.2 cm 
\noindent
{\bf Remark:} Each pole  $\xi_k(j)$  lies in a disk centered in $x_k(j)$ of radius ${1\over 2} [x_k(j)]^2+O( [x_k(j)]^3)$, as depicted in figure \ref{poles}. The zeros are ordered as $|x_{k+1}(j)|<|x_k(l)|<|x_k(j)|$, where $(j,l)=(1,2)$ or $(2,1)$. Consider disks centered at these zeros of radius $|x_{k}(j)|^2$, $|x_k(l)|^2$, 
$|x_{k+1}(j)|^2$. $K$ is constructed in the proof of the theorem  in such a way that when  $k\geq K$ the disks do not intersect. Thus, Theorem \ref{theo1} is consistent, because  $\xi_k(j)$ is closer to $x_k(j)$ than to any other $x_{k^\prime}(j^\prime)$. 

\begin{figure}
\epsfxsize=10cm

\epsfysize=20cm

\centerline{\epsffile{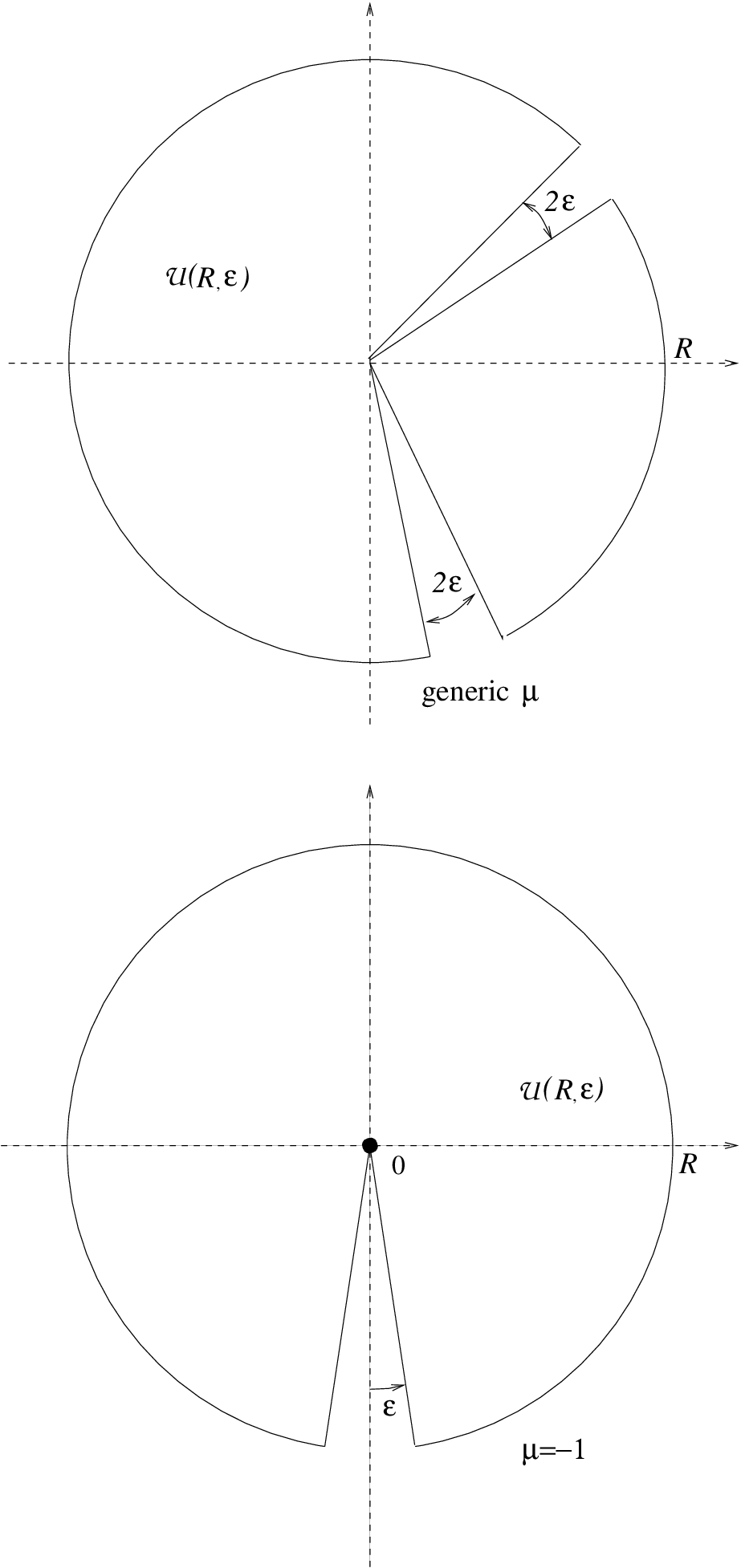}}
\caption{}
\label{disk}
\end{figure}

\vskip 0.3 cm
 We may ask if there are other  poles of $y(x)$, in the domain of convergence of (\ref{pro2y}), other than 
 those    of Theorem \ref{theo1}.
 Let  $R<$ radius of convergence of (\ref{pro2y}) $<1$ be the radial coordinate for a point on the ray where
 the zeros of $y_1(x)$ lie.   We prove that  {\it $y(x)$ has no poles that are more than $R^2$-distant from
 the rays}, as in figure \ref{poles}.  To formalize this statement, let  ${\cal U}(R,\epsilon)$ be a close 
 domain constructed by taking a  disk centered at $x=0$ of radius $R$, minus two sectors bisected by the
 rays where the zeros lie.  Each sector has  angular amplitude  $2\epsilon$. See
 figure \ref{disk}, where the general case and the case $\mu=-1$ are depicted. 

% % % % % % % % % % % % % % % %% % % %
%%
%
%
%%
% % % % % % % % % % % % % % % %

\bth
\label{theo3}
Let $y(x)$ be (\ref{pro2y}). For any  small $\epsilon>0$ , there exist $R_\epsilon<$ radius of convergence of (\ref{pro2y}), such that   $y(x)$ has no poles in ${\cal U}(R_\epsilon, \epsilon)$. When $\epsilon\to 0$,    
$R_\epsilon$ can be chosen to be any number such that:  
$$
R_\epsilon< {|2\mu-1|\over C_f}~\tan \epsilon~(1+O(\tan^2\epsilon)),~~~~~\mu\neq {1\over 2}.
$$ 
$$
R_\epsilon< {2\nu^2\over C_f} \tan^2(\epsilon)~ (1+O(\tan^3 (\epsilon))), ~~~~~\mu = {1\over 2}.
$$
\eth

\qed

\vskip 0.2 cm
  
 The asymptotic estimate  $R_\epsilon\sim {|2\mu-1|\tan \epsilon\over C_f}$ for $\epsilon \to 0$  means that the poles of $y(x)$, if they exist, get closer to the the rays where the zeros  lie, as their absolute value decreases. Their distance from the the rays is at most 
${|2\mu-1|\epsilon\tan \epsilon/ C_f}\cong {|2\mu-1|\epsilon^2/ C_f}\cong {C_f R_\epsilon^2/ |2\mu-1|}$ (or $\sqrt{C_f}R_\epsilon^{3\over 2}/\sqrt{2}\nu$ if $\mu={1\over 2}$).

 \vskip 0.2 cm 
 A relevant example of Frobenius manifold is given by the {\it quantum cohomology of the two dimensional complex projective space} $CP^2$. A solution of $PVI_\mu$, with $\mu=-1$, is associated to this manifold  \cite{Dub}.  It is a solution of the form (\ref{pro2y}) (please, see section \ref{PVIFROBENIUS} in  order to understand this fact).  In this case, we prove that 
$$
\nu= {2\ln {\bf G}\over \pi},~~~~{\bf G}={1+\sqrt{5}\over 2}
$$
where ${\bf G}$ is the {\it golden ratio}.  This very special value  makes the two rays of zeros coincide with the negative immaginary axis. Note that we can choose $\vartheta=\pi$, in such a way that (\ref{pro2y}) becomes a branch with branch cut along the negative axis. The result is summarized in the following theorem: 

%%%%%%%%%%%%%%%%%%%%%%  THEOREM %%%%%%%%%%%%%%%%%%%%%%%%%%%%%%%%%
\bth
\label{theo2}
The branch (\ref{pro2y}) associated to the Quantum Cohomology of $CP^2$,  defined for $|\arg x|<\pi$,  which satisfies the equation $PVI_{\mu=-1}$, is identified by the following integration contants: 
$$
\nu = {2 \ln {\bf G}\over \pi}=0.30634...
$$

\be
\label{3vetrini}
d={i\over 2} \ln
\left\{-
\pi^2~{({\bf G}^4+1)^2\over ({\bf G}^2+1)^2}~2^{16i\nu}~{(1-2i\nu)^2 ~\nu^2 \over (1+2i\nu)^2}~{\Gamma(1-2i\nu)^4\over 
\Gamma(1-i\nu)^8}
\right\},
\ee
The first approximation $y_1(x)$ has  two  infinite sequences of zeros  accumulating at $x=0$ along the negative imaginary axis:   
\be
\label{synzeros}
  x_k(j)= -i\exp
\left\{
-{\Re d\over \nu} -{2(j-1)\over \nu}\left|\arccos{3\over \sqrt{4\nu^2+9}}\right|\right\}\exp\left\{-{k\pi\over \nu}\right\}, ~~k\in{\bf N}, ~~~j=1,2.
\ee
The branch (\ref{pro2y}) 
 has two sequences of poles $\{\xi_k(1)\}_{k\geq k_0}$, 
$\{\xi_k(2)\}_{k\geq k_0}$ in a neighborhood of $x=0$, which accumulate at $x=0$ as $k\to \infty$, asymptotically approaching the negative imaginary axis, according to the asymptotic expansion (\ref{POLESintro}), 
where the first terms are: 
$$
\Delta_3(1)={176\nu^4+185+352\nu^2 \over 1024 (\nu^2+1)^2}=0.1792...,~~~
\Delta_3(2)={401+176\nu^4+352\nu^2\over 1024(\nu^2+1)^2}=0.3555...
$$
$$
\Delta_4(1)=-{57+48\nu^4+96\nu^2\over 1024(\nu^2+1)^2}=-0.05422...,~~~
\Delta_4(2)=-{273+48\nu^4+96\nu^2)\over 1024(\nu^2+1)^2}=-0.2305...
$$
\eth
\qed

%%%%%%%%%%%%%%%%%%%% end theorem %%%%%%%%%%%%%%%%%%%%%%%%%%%%%%%%%%%%%%%%%%%%

\vskip 0.3 cm 
\noindent
{\bf Corollary:} {\it
 $d$ of Theorem \ref{theo2} has the following series 
$$
d=
{\pi\over 2}-8\nu\ln(2)+2\arccos{1\over \sqrt{1+4\nu^2}}+4 \sum_{n=1}^\infty {(-1)^n(1-4^n)\zeta(2 n+1)\over 2n+1}  ~ \nu^{2 n+1}~
+~i~{\pi\nu\over 2} +k\pi
$$
where  $
2\arccos\bigl({1/ \sqrt{1+4\nu^2}}\bigr)=2\sum_{n=0}^\infty (-1)^n(2\nu)^{2n+1}
$, and $k\in{\bf Z}$. 
The above series is absolutely convergent for $|\nu|<{1\over 2}$, therefore well defined for $\nu = {2 \over \pi}\ln {\bf G}$.  In particular, $
\Im d={\pi\nu\over 2}
$. 
}
\qed
\vskip 0.3 cm 

According to (\ref{Nishigamo}), we take $k=0$ in the corollary. The first terms of the series then give the numerical value $
  \Re d= 1.08323...$ and  ${\Re d\over \nu}= 3.53595...
$. 
Other  numerical values are: $\arg\left({3-2i\nu\over 3+2i\nu}\right)= -2 ~\hbox{arcos}\left({3\over \sqrt{4\nu^2+9}}\right)=-0.402...$, therefore  $e^{{1\over \nu} \arg\left({3-2i\nu\over 3+2i\nu}\right)}=0.268...$, $e^{-{\Re d\over \nu}}=0.0291...$, $e^{-{\pi\over \nu}}=3.52...\times 10^{-5}$. It follows that the first two zeros for $k\geq 1$  are  $x_0(1)=-2.91...\times 10^{-2}~i$ and $x_0(2)=-7.82...\times 10^{-3}~i$.

\vskip 0.2 cm 
 In Section \ref{PVIFROBENIUS}, we review the dependence of the integration constants $d$ and $\nu$ on the 
monodromy data associated to a solution of $PVI_\mu$, for any $\mu$.  It turns out that $d$ depends 
explicitly on $\nu$.
 In Section \ref{logLIM} we expand $d$ as a convergent Taylor series of $\nu$:
$$
  d= \sum_{n=1}^\infty d_n \nu^n,~~~~~d_n\in{\bf C}
$$
  This implies that the zeros of $y_1(x)$ shrink to $x=0$, when $\nu\to 0$. Namely, for $j=1,2$ and $k\geq 0$, we have
$$
 x_k(j)\to 0~~~\hbox{ when }\nu\to 0
$$
We prove that,  for $x\neq 0$, the following holds: 
\be
\label{Pilog}
 \lim_{\nu\to 0} y_1(x)=\left({1\over 2}-\mu\right) \left[\left(\mu-{1\over 2}\right)(\ln(x)+d_1)-2\right](\ln x+d_1) 
\ee
In the same way,   $y_n(x)$ converges 
to a polynomial of $\ln x$, when $\nu \to 0$. We have proved this up to $n=4$ and conjecture that this is true for any $n$. Therefore, we expect that, when $x\neq 0$,  the limit of $y(x)$ for $\nu \to 0$ exists, whith asymptotic expansion 
\be
\label{conjintro}
\lim_{\nu \to 0} y(x)\sim\left[\sum_{n=1}^\infty x^{n-1} P_n(\ln x)\right]^{-1},~~~x\to 0 
\ee
where $P_n(\ln x)$ are certain polynomials of $\ln x$. 
On the other hand, we showed in \cite{D1} that, for $\nu=0$,  $PVI_{\mu}$ has  solutions with asymptotic expansion coinciding with the right hand side of the above (\ref{conjintro}). This verifies the conjecture.  The Chazy solutions of \cite{MazzChazy} are also  re obtained from the limit of solutions (\ref{pro2y}) for $\nu \to 0$. 
 
In Section \ref{exempicard} the example of Picard solutions is didactically discussed, and we prove convergence of (\ref{POLESintro}). In Section \ref{generalPVI} we sketch the case of the general PVI equation. The rest of the paper is devoted to the proof of the theorems.

%  % % % % % % % % % % % %% %   

\section{Parameterization in terms of Monodromy Data}
\label{PVIFROBENIUS}

The equation $PVI_{\mu}$ is associated to a semisimple Frobenius manifold of complex dimension 3, locally described by a system  of canonical complex  coordinates $u_1$, $u_2$, $u_3$. It is known that a branch of a solution is parameterized by the monodromy data of the manifold, namely $\mu$ and the three entries $s_{12}$, $s_{13}$, $s_{23}$  of the {\it Stokes' matrix} of the manifold, while the independent variable is $x=(u_3-u_1)/(u_2-u_1)$. These results are estabilished in  \cite{Dub}. 
A remarkable example is  the Frobenius manifold given by the  {\it Quantum Cohomology of $CP^2$} \cite{Dub}. In this case, 
$$
\mu=-1,~~~~~s_{12}=s_{13}=s_{23}=3
$$
 (for the compution of the stokes matrix of $CP^d$ see  \cite{Dub} when $d=2$,  and  \cite{guzz} when $d\geq 3$).

\vskip 0.2 cm 
In place of $\mu$ and  the  entries of the Stokes matrix, we will use the equivalent quantities 
$$
   \theta_\infty:=2\mu,~~~~~p_{0x}:=2-s_{13}^3,~~~p_{01}:=2-s_{12}^2,~~~p_{x1}:=2-s_{23}^2
$$
The above are  usually employed in the isomonodromy preserving deformation approach to PVI, estabilished  in \cite{JMU} and \cite{Jimbo}.  They are monodromy data of the following $2\times 2$ Fuchsian system of ODE associated to $PVI_\mu$:
\be
   {d\Psi\over d\lambda}=\left[ {A_0(x)\over \lambda}+{A_x(x) \over \lambda-x}+{A_1(x)
\over
\lambda-1}\right]\Psi,~~~\lambda\in{\bf C}.
\label{SYSTEM}
\ee
$$
 A_0+A_1+A_x = \hbox{diag}\left(-{\theta_{\infty}\over 2},{\theta_{\infty}\over 2}\right),~~~
\hbox{ Eigenvalues}~( A_i) =0, ~~~i=0,1,x;
$$
A ordered base of loops $\Gamma$ is chosen in the fundamental group of ${\bf C}\backslash \{0,x,1\}$. The three basic loops encircle  $0,x$ and $1$ respectively.   The matrices $A_i(x)$ depend on $x$ in such a way that the monodromy  group w.r.t.  $\Gamma$ is independent of small deformations of  $x$.  By small deformation it is meant that $x$ does not go around a loop around 0 or 1, namely that some banch cuts are chosen. For example, we may choose $-\pi \leq \arg x <\pi$ and $-\pi \leq \arg (x-1) <\pi$. Correspondingly, the solution $y(x)$ of $PVI_\mu$ is to be regarded as {\it a branch}. For the local analysis around $x=0$, we just consider $-\pi \leq \arg x <\pi$. 
 Let $M_0,M_x,M_1$ be the monodromy matrices of a fundamental solution of the Fuchsian system w.r.t $\Gamma$. Let $M_\infty=(M_1M_xM_0)^{-1}$ be the monodromy at infinity. According to \cite{JMU} and \cite{Jimbo} the quantities $p_{ij}$ are: 
 $$p_{ij}=\hbox{Tr} (M_iM_j),~~~j=0,x,1.
$$ 
Let also 
$$
p_\infty:= \hbox{Tr} M_\infty= 2\cos(\pi\theta_\infty),~~~~~
$$
 The $p_{ij}$'s and $p_\infty$  are  coordinates for  the space of monodromy  data of the  class of Fuchsian systems above. This space  is an affine cubic surface \cite{Iwa} \cite{Jimbo}:
\be
\label{cubic}
p_{0x}^2+p_{01}^2+p_{x1}^2+p_{0x}p_{01}p_{x1}-2(p_{0x}+p_{01}+p_{x1})(2+p_\infty)+8+p_\infty^2+8p_\infty=0
\ee
Note that $\theta_\infty=2\mu$ is fixed by the equation and
 $p_{0x},p_{01},p_{x1}$ are not independent, because of the cubic relation. Accordingly, only two
 complex parameters are free.  The following  facts follow from the general theory of Painlev\'e VI (see \cite{guz2010}, section 2 \footnote{
%
% FOOTNOTE
%
Keep into account that,  for $PVI_\mu$, each of the monodromy matrices  $M_0,M_x$ and $M_1$ has a Jordan form $\pmatrix{1 & 2\pi i\cr 0 & 1}$ (but they cannot be put simultaneously in upper triangular form, in general), so they are not the identity matrix $I$,
}
): 

\vskip 0.2 cm 
\noindent
{\it Let the basis of loops $\Gamma$ be fixed. 
 
1) If the monodromy group $<M_0,M_x,M_1>$ is not reducible, or  $M_\infty\neq I$, the above $p_\infty$, $p_{ij}$'s are a good system of coordinates for the monodromy group \cite{Iwa},\cite{Jimbo}. 

2) If $M_\infty\neq I$, there is a one to one correspondence between a branch $y(x)$ and a point in the space of monodromy data. 

3) If the monodromy group $<M_0,M_x,M_1>$ is not reducible and  $M_\infty\neq I$, a branch $y(x)$ of a transcendent of $PVI_\mu$ is uniquely parameterized by
the $p_\infty$'s (i.e. $\theta_\infty$)  and $p_{ij}$'s, to which it is in one to one correspondence.
}
\footnote{
When the monodromy group is reducible, but  $M_\infty$ is not the identity, the one to one correspondence between a point in the space of monodromy data and a branch still holds, but the $p_\infty$, $p_{ij}$'s are no longer a good parametrization. The solutions in this case are known (see the Riccati solutions \cite{watanabe}, \cite{mazzoccoratio}). 
}

\vskip 0.3 cm
 As a consequence of 3) above, the two complex integration constants of the branch $y(x)$  are functions of $p_\infty$
 (i.e. $\theta_\infty$) and $p_{ij}$ . 
A remarkable fact, established in  Jimbo's paper \cite{Jimbo}, is that {\it this parametrization is explicit}, namely the integration constants are elementary or classical  transcendental functions of $p_\infty$, $p_{ij}$, $i,j=0,x,1$.  Jimbo computed the parametrization for the generic PVI and most of the range of $p_{ij}$, except for  $p_{ij}<-2$.  \footnote{As a consequence of this explicit parameterization of the three couples of
 integration constants at the three critical points in terms  of {\it the same} monodromy data, the
 connection problem is solved.   This is precisely the power of the method of
 monodromy preserving deformations.} 

\vskip 0.2 cm 

 Solutions (\ref{pro2y}) occur when $s_{13}>2$, namely $p_{0x}<-2$. In this case, the parametrization of $\nu$ and $d$ in terms of the monodromy data is computed in \cite{guz2010}, and it is summarized in the following:

%%%%%%%%%%%%%%%%%%%%%%%% PROPOSITION %%%%%%%%%%%%%%%%%%%%%%%%%%%%%%%%%%%%%%%%
\bpr
\label{lem1}
Let $PVI_\mu$ be given, namely let $\mu$ be given. 
The branch  (\ref{pro2y})  is associated to monodromy data such that $p_{0x}<-2$. 
  The  real integration constant  $\nu>0$ is obtained from
$$
 p_{0x}=-2\cosh(2\pi\nu)~~~\hbox{(namely, $
  s_{13}=2\cosh(\pi \nu) $)}
$$
Therefore: 
\be
\label{nugen}
\nu= {2\ln G\over \pi},~~~~~G:=\left[{\sqrt{p_{0x}^2-4}-p_{0x}\over 2}\right]^{1\over 4}>1,
\ee
$\diamond $ If $2\mu\neq 2i\nu+2m+1$, $m\in{\bf Z}$, the complex integration constant $d$ is:
\be
\label{generald}
d= {i\over 2} \ln
 \left\{-
{4~16^{2i\nu}~\Gamma({3\over 2}-{\mu}-i\nu)^2 ~\Gamma({\mu}+{1\over 2}-i\nu)^2 
\over \Bigl(2\nu+i(1-2\mu)\Bigr)^2\nu^2~\sinh(2\pi\nu)^2~\Gamma(-i\nu)^4}~\times \right.~~~~~~~~~~~~~~~~~~~~~~~~~~~~~~~~
\ee
$$~~~~~~~~~~~~~~~~~~~~
\times 
\left.
\left[~{1\over 2}
\Bigl(e^{2\pi\nu}p_{x1}-p_{01}\Bigr)\sinh(2\pi\nu)+\Bigl(\cos(2\pi\mu)+1\Bigr)\Bigl(e^{2\pi\nu}+1\Bigr)~
\right]~
\right\}.
$$

\noindent
$\diamond$ If $2\mu= 2i\nu+2m+1$, $m\in{\bf Z}$, then $p_{0x}=2\cos(2\pi\mu)$, $p_{01}-2=(p_{x1}-2)e^{\pm2\pi\nu}$, and $d$ is as below:

\vskip 0.2 cm 
-- If  $2\mu=2i\nu+1-2m$, $m=1,2,3,...$: 
\be
\label{generald1}
d=-{i\over 2}\ln
\left\{{\nu^2\Gamma(m)^2\Gamma(2i\nu+1-m)^2
\over
16^{2i\nu}\Gamma(1+i\nu)^4
}(p_{x1}-2)
\right\}
\ee
$$
2-p_{01}=(2-p_{x1})e^{2\pi\nu}
$$
\vskip 0.2 cm 
-- If  $2\mu=2i\nu+1-2m$, $m=-1,-2,-3,...$: 
\be
\label{generald2}
d={i\over 2}
\ln\left\{
{
\sinh(\pi\nu)^4~16^{2i\nu}\Gamma(1-m)^2\Gamma(1+i\nu)^4\Gamma(m-2i\nu)^2
\over
\nu^2\pi^4
}(p_{x1}-2)
\right\}
\ee
$$
2-p_{01}=(2-p_{x1})e^{-2\pi\nu}
$$
\vskip 0.2 cm 
-- If  $2\mu=-2i\nu+2m-1$, $m=1,2,3,...$ (but $m\neq 1$ in (\ref{pro2y})):
\be
\label{generald3}
d=-{i\over 2}
\ln\left\{
{
\nu^2\Gamma(2i\nu+1-m)^2\Gamma(m)^2
\over 
16^{2i\nu}\Gamma(1+i\nu)^4
}(p_{x1}-2)
\right\}
\ee
$$
2-p_{01}=(2-p_{x1})e^{2\pi\nu}
$$
\vskip 0.2 cm 
-- If  $2\mu=-2i\nu+2m-1$, $m=0,-1,-2,-3,...$:
\be
\label{generald4}
d={i\over 2}
\ln\left\{
{\sinh(\pi\nu)^4~16^{2i\nu}\Gamma(1-m)^2\Gamma(1+i\nu)^4\Gamma(m-2i\nu)^2
\over
\nu^2\pi^4
}(p_{x1}-2)
\right\}
\ee
$$
2-p_{01}=(2-p_{x1})e^{-2\pi\nu}
$$
Being $\nu\neq 0$, in the above cases $p_{01}\neq p_{x1}$. 
\epr

%%%%%%%%%%%%%%%%%%%%%%%%%%%%% END PROPOSITION %%%%%%%%%%%%%%%%%%%%%%%%%%%%%%%%%%%%%%%%%%%%

\noindent
{\it Proof:}  Section \ref{monodromydata}.

\vskip 0.3 cm 
Note that the freedom in the choice of the branch of the logarithm defines $d$ up to  $d\mapsto d+k\pi$, $k\in {\bf Z}$. Such freedom does not affect $y(x)$ and we can choose $
0\leq \Re d\leq \pi$, as in (\ref{Nishigamo}).

 For a given $PVI_\mu$, let $\nu>0$ and $d\in{\bf C}$ be given (equivalently, let monodromy data be given). Let us denote  
$$
y(x,\nu,d):=y(x) \hbox{ with behavior  {\rm (\ref{pro2y})} with }
-\pi\leq \arg x <\pi
$$ 
This is the branch with behavior (\ref{pro2y}). Its {\it analytic continuation} when $x$ goes around a small loop around $x=0$, or $x=1$, or $x=\infty$ may have a behavior different from (\ref{pro2y}).  For the local analysis  at $x=0$, it is enough to consider the analytic continuation when  $x$  goes around a loop around $x=0$, namely $x\mapsto x e^{2\pi i}$ ($|x|<1$).  The new branch  is parametrized in terms of new monodromy data $p_{ij}^\prime$, $i,j\in\{0,x,1\}$, computable by an action of the braid group  as follows (see \cite{DM}): 
$$
  p_{0x}^\prime=p_{0x},~~~~~p_{01}^\prime=
-p_{01}+2p_{x1}\cosh(2\pi\nu)+4\cos(2\pi\mu)+4
$$
$$
 p_{x1}^\prime=p_{x1}(4\cosh(2\pi\nu)^2-1)-2\cosh(2\pi\nu)p_{01}+4(\cos(2\pi\mu)+1)(2\cosh(2\pi\nu)+1)
$$
We see that $p_{0x}^\prime=p_{0x}<-2$, thus the new branch has again a behavior (\ref{pro2y}) and Proposition \ref{lem1} holds. In particular, $\nu$ is unchanged.  As a consequence, we have: 

\bpr
 The analytic continuation of the branch $y(x,\nu,d)$ corresponding to the loop $x\mapsto xe^{2\pi i}$ is:
\be
\label{analcont}
y(x,\nu,d)\mapsto   y(x,\nu,d+2\pi i\nu)=y(xe^{2\pi i},\nu,d),~~~-\pi\leq \arg x <\pi.
\ee
\epr

\noindent
{\it Proof:} 
Substituting into the formulas of Proposition \ref{lem1} the expressions of $p_{x1}^\prime$ and $p_{01}^\prime$  and simplifying,  we find that the new $d$, denoted $d^\prime$, is $
 d^\prime= d+2\pi i\nu.
$. 
Observe also that that $\exp\{2id^\prime\}x^{2i\nu}=\exp\{2id\}(xe^{2\pi i})^{2i\nu}$.  Immediately (\ref{analcont}) follows. 
\qed

\vskip 0.3 cm 
\noindent
{\bf Remark on the poles of a branch:} The rays of Theorem \ref{theo1}, where the zeros  $x_k(1)$ and $x_k(2)$ lie,  
  may  be outside the range $-\pi\leq \arg x<\pi$,  depending on their angles 
${\Im d\over \nu}$ and ${\Im d\over \nu}+{1\over \nu} \ln \left| {2\mu-1+2i\nu\over 2\mu-1-2i\nu}\right|$ (note: the expansion (\ref{pro2y}) is defined for $|\arg x|<\vartheta$, 
with  $\vartheta>0$ and $|x|\neq 0$ sufficently small. Recall that $\vartheta$ is arbitrary (but fixed) and this fact allows to find zeros with angles outside $-\pi\leq \arg x<\pi$). When this happens,  one or both the sequences of  the zeros do not fall in the domain  $-\pi\leq \arg x<\pi$.   Accordingly, the branch  $y(x,\nu,d)$ does not have poles. 
The analytic continuation of the branch when $x$ goes around a loop $x\mapsto xe^{2\pi i}$ ($|x|<1$) is (\ref{analcont}). 
The shift $d \mapsto d+2\pi i\nu$ changes the immaginary exponent of the $x_k(1)$'s by $-2\pi i$. This implies that, by a sufficient number of loops,  we can always find a branch with poles, namely such that   at least one  of the two sequences of zeros is  in the range $-\pi \leq \arg x<\pi$.

\vskip 0.2 cm

% % % % % % % % %  % % %

\section{Limit for $\nu\to 0$}
\label{logLIM}
 Suppose that $d$ vanishes with $\nu$ as $\nu\to 0$. Namely: 
\be
\label{limitofd}
  d=d_1\nu\Bigl(1+d_2(\nu)\Bigr),~~~d_1\in{\bf C}\backslash\{0\},~~~d_2(\nu)\to 0,~~~\nu\to 0.
\ee
 If this happens,  the zeros $x_k(j)$ of Theorem \ref{theo1}  shrink to $x=0$, provided that $k\geq 0$. \footnote{
As for $x_k(2)$, use the fact that: 
$$
\exp\left\{
-{i\over 2}\ln
\left(
{2\mu-1+2i\nu\over 2\mu-1-2i\nu}
\right)
\right\}=\exp\left\{
-{i\over 2}\ln
\left(1+{4\nu i\over 2\mu-1}+O(\nu^2)
\right)
\right\}={4\over 2\mu-1}+O(\nu),~~~\nu\to 0.
$$
}

\bpr
\label{seki}
Suppose that $d$ is as in (\ref{limitofd}).  
Then, for $x\neq 0$, there exists: 
$$
\lim_{\nu\to 0} y_1(x)= P_1(\ln x)
$$
where $P_1(\ln x)$ is the polynomial of $\ln x$ appearing in the r.h.s. of (\ref{Pilog}).
\epr

\noindent
{\it Proof:} Substitute $d=d_1\nu(1+d_2(\nu))$ in (\ref{firstordermu}). Then expand $
  \exp\{2id\}=1+2id_1\nu+o(\nu)$ and $x^{2i\nu}=1+2i\nu\ln x+O(\nu^2)$. 
The structure of the coefficients  $A_{1,-1}$, $A_{10}$, $A_{11}$ allows simplification of the divergences $\nu^{-2}$ and $\nu^{-1}$ contained in the coefficients themselves. Therefore, $y_1(x)$ is expanded in series for $\nu\to 0$, which by direct computation is easily verified to be:
$$
  y_1(x)= P_1(\ln x)+\nu(\hbox{\rm polynomial $\ln x$})+\nu^2(\hbox{\rm polynomial $\ln x$})+...,~~~~~x\neq 0
$$
\rightline{\qed}

\vskip 0.2 cm 
\noindent
{\bf Remark:} If  $
 d=d_0+d_1\nu+d_2\nu^2+...$, $d_0\neq 0$, 
then 
$$
y_1(x)={(2\mu-1)^2(\cos(2d_0)-1)\over 8 \nu^2}-{(2\mu-1)\left((2\mu-1)(\ln x+d_1)-2\right)\sin(2d_0)\over 4\nu}+O(\nu^0),~~~\nu\to 0
$$
Thus, the limit for $\nu\to 0$ and $x\neq 0$ does not exist. 

\vskip 0.2 cm 
The same computation of the proof of Proposition \ref{seki}  can be done for $y_n(x)$. We verified, up to $n=4$, that the divergences in $A_{nm}$ for $\nu\to 0$ are canceled by the expansion of $\exp\{2i(d_1\nu(1+d_2(\nu))\}$ and  $x^{2i\nu}=\exp\{2i\nu\ln x\}$ for $\nu\to 0$, so that there exists:
$$
\lim_{\nu\to 0} y_n(x)= P_n(\ln x)
$$
  where:
$$
P_n(\ln x)= (2\mu-1) p_0(\mu,d_1)+\sum_{N=1}^{2n}(2\mu-1)^Np_N(\mu,d_1)~\ln^N(x)
$$
The $p_N(\mu,d_1)$ are certain polynomials of $\mu$ and $d_1$. For example: 
$$
P_2(\ln x)= -{(2\mu-1)^4\over 32}\ln^4(x)-{(2\mu-1)^3((2\mu-1)d_1-2)\over 8}\ln^3(x)+
$$
$$
-3\left(\mu-{1\over 2}\right)^2\left({5\over6}+\left(\mu-{1\over 2}\right)^2d_1^2+(1-2\mu)d_1\right)\ln^2(x)+
$$
$$
-(2\mu-1)\left(\left(\mu-{1\over 2}\right)^3d_1^3-3\left(\mu-{1\over 2}\right)^2d_1^2+{5\over 4}\left((2\mu-1\right)d_1+{1\over2}\mu-{3\over 4}\right)\ln(x)+
$$
$$ -(2\mu-1)\left(2+\left(\mu-{1\over 2}\right)^2d_1^3+{3\over 2}\left(1-2\mu\right)d_1^2+2d_1\right)\left(\left(\mu-{1\over 2}\right)d_1-1\right)
$$
Vanishing of $d$ as $\nu \to 0$ is not an arbitrary assumption, because the following holds:

\bpr
\label{limitd}
First case: suppose that $\mu\not \in \{ \pm 2i\nu+2m+1, ~2m+1\}_{m\in {\bf Z}}$. The integration constant $d$, given by (\ref{generald}), is expanded as a Taylor series, convergent for $\nu$ sufficiently small:
$$
d=\sum_{n=1}^\infty d_n\nu^n,~~~~~\nu\to 0,~~~d_n\in{\bf C}
$$
where
$$
d_1= 2\gamma-4\ln2+{i\pi\over 2}+2~ \psi\left(\mu+{1\over 2}\right)-\pi \tan(\pi \mu)+{i\pi(p_{x1}-p_{01})\over 8 \cos^2(\pi\mu)}
$$

Second case: when  $\mu=\pm 2i\nu+2m+1$, $m\in {\bf Z}$. The integration constant $d$ is expanded as a Taylor series, convergent for $\nu$ sufficiently small, as follows:
\vskip 0.2 cm 
\noindent
 -- The  cases of (\ref{generald1}) and  (\ref{generald3}): 
$$
d=-{i\over 2}\ln\left({2-p_{x1}\over 4}\right)+2\left(\psi(m)+\gamma-2\ln 2\right)\nu+\sum_{n=2}^\infty d_n\nu^n,
$$
$$ 
m=1,2,3,...,~~~~~~~d_n\in{\bf C}
$$
\vskip 0.2 cm 
\noindent
-- The  case of (\ref{generald2}) and  (\ref{generald4}): 
$$
d={i\over 2}\ln\left({2-p_{x1}\over 4}\right)+2\left(\psi(1-m)+\gamma-2\ln 2\right)\nu+\sum_{n=2}^\infty d_n\nu^n,
$$
$$ 
m=0,-1,-2,-3,...,~~~~~~~d_n\in{\bf C}
$$
In the above formulas $
\psi(m)+\gamma= \sum_{k=1}^{m-1}{1\over k}
$. 
\epr 

\noindent
{\it Proof:} Section \ref{prooflimitd}. 

\vskip 0.2 cm 

The above results allow us to formulate the following:

\vskip 0.3 cm 
\noindent
{\bf Conjecture:} {\it 
Let $\mu\neq {1\over 2}$ and assume that  $d$ is as in (\ref{limitofd}). 
Let $y(x)$ be  the branch of (\ref{pro2y}). Then, for $x\neq 0$, there exist the limit of $y(x)$ for $\nu \to 0$, with  asymptotic series  (\ref{conjintro}). 
}

\vskip 0.3 cm 
Recall that $p_{0x}=-2\cosh(2\pi\nu)$. If $\nu\to 0$, then $p_{0x}\to -2$. 
The result established in \cite{D1} (see also \cite{guz2010} for a review) states that  to the monodromy data such  $p_{0x}=-2$ and $p_{x1}$, $p_{01}$ arbitrary, there is associated a branch of a solution of $PVI_{\mu}$, $\mu \neq {1\over 2}$, with critical behavior:
\be
\label{logbehavior}
  y(x)={1\over P_1(\ln x)}\left(1+O(x\ln^2(x))\right)= -{4\over (2\mu-1)^2\ln^2 (x)}\left(1+O\left({1\over \ln x}\right)\right),~~~x\to 0
\ee
The full asymptotic expansion can be computed (by substitution into the equation) and it coincides with  (\ref{conjintro}).  This shows that the conjecture is true. 
 Note that the monodromy data, to which a branch is in one to one correspondence, are contained in $d_1$,
 when $\mu\neq \pm 2i\nu+2m+1$. Also note that the first term in the critical behavior (\ref{logbehavior})
 does not depend on $d_1$, namely on the monodromy data $p_{x1}$, $p_{01}$.

\begin{figure}
\epsfxsize=10cm

\epsfysize=10cm

\centerline{\epsffile{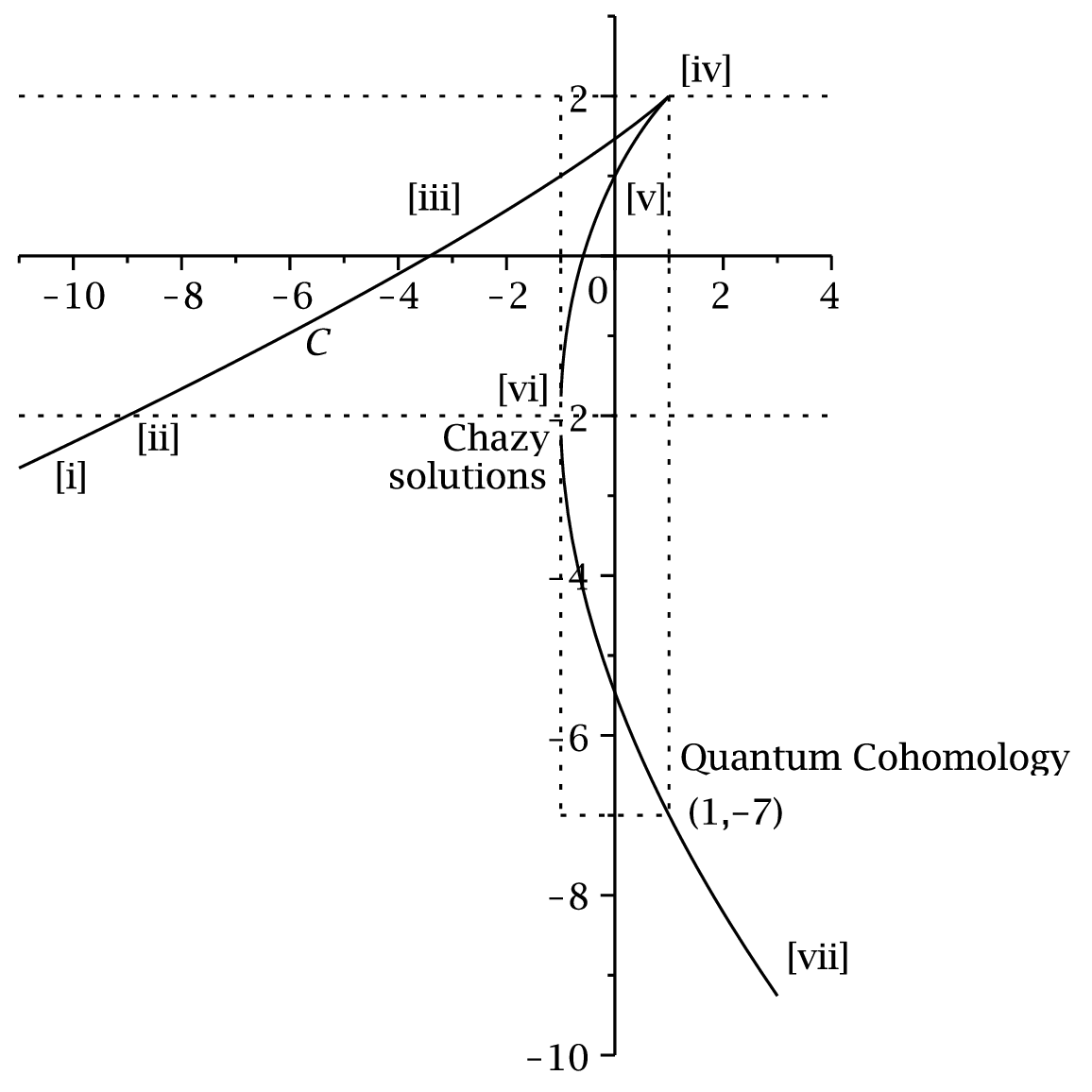}}
\caption{Cubic curve in the plane $(\cos\pi\theta_\infty,p_{0x})=(\cos(2\pi\mu),p_{0x})$}
\label{CUBIC}
\end{figure}
%
%%

% % %%  %%  %    %% % %  % %% % % % % % % % % % % % % %   %% % % % %

\subsection{Chazy Solutions, $\nu \to 0$ and  $\mu\to-{1\over 2}$}

In the special cases $\mu=\pm 2i\nu+2m+1$, $d$ does not vanish, therefore  the limit of $y_1(x)$ for $\nu \to 0$ and $x\neq 0$ diverges as $\nu^{-2}$, except possibly when $p_{x1}=-2$.   
Also note that the restriction $\mu\neq {1\over 2}+m$, $m\in{\bf Z}$ avoids divergences of $d$  in Proposition \ref{limitd}.  Such cases may be analyzed separately. Preliminarily, we observe that 
the symmetries of PVI imply that one can always assume  that $-1\leq \Re \mu <0$. 
%\be
%\label{restricted}
 %  -1\leq \Re \mu <0
%\ee
Thus, it is enough to analyze the case $\mu=-{1\over 2}$. We will limit ourselves to the analysis of  the solutions  associated to monodromy data which satisfy the condition: 
$$
p_{x1}=p_{01}=p_{0x}=-2\cosh(2\pi\nu)
$$
The cubic relation (\ref{cubic}) becomes: 
\be
\label{cubic1}
3p_{0x}^2+p_{0x}^3-12~p_{0x}\Bigl(1+\cos(2\pi \mu)\Bigr)+4\cos^2(2\pi\mu)+16\cos(2\pi\mu)+8=0
\ee
This curve is represented in figure \ref{CUBIC} by the continuous line, when $\cos(2\pi\mu)$ is real. Observe that $\nu\geq 0$ is well defined when $p_{0x}\leq -2$, This condition singles out two portions of the curve, namely [i]$\cup$[ii] and [vi]$\cup$[vii]. If  we require that also $\mu$ is real, namely $-1\leq \cos(2\pi\mu)\leq 1$, the condition singles out only the portion  of the curve between $(-1,-2)$ and $(1,-7)$. In this case, the functional relation between $\mu$ and $\nu$  is established in Lemma \ref{tetanu} of Section 
\ref{proofchazyd}: 
$$
   \mu= -{1\over 2} +{1\over 2\pi}~ \hbox{\rm arg}\left(c_1(\nu)-{i\over 2}\sqrt{c_2(\nu)}\right),
$$

\be
\label{starrr}
\left\{
\matrix{
 c_1(\nu)&:=& 2+3\cosh(2\pi\nu)-\cosh(3\pi\nu)-3\cosh(\pi\nu),
\cr
\cr
c_2(\nu)&:=& 96\cosh(\pi\nu)+52\cosh(3\pi\nu)+12\cosh(5\pi\nu)+
\cr
& & -50-78\cosh(2\pi\nu)-30\cosh(4\pi\nu)-2\cosh(6\pi\nu)
}
\right.
\ee 

\noindent
where $
0\leq \nu \leq {2  \ln {\bf G} / \pi}$ and ${\bf G}=(1+\sqrt{5})/ 2$.  The discriminat $c_2(\nu)\geq 0$ and the square root is the positive one. 
When we choose the argument with determination
$$
-\pi\leq \hbox{\rm arg} \left(c_1(\nu)-{i\over 2}\sqrt{c_2(\nu)}\right)\leq 0
$$
\noindent
it follows that $\mu$ can be expanded as a convergent series, as in the following proposition. 

\bpr 
\label{chazyd}
Let $
p_{x1}=p_{01}=p_{0x}=-2\cosh(2\pi\nu)\leq -2
$ and let $\mu$ be real, $-1\leq \mu<0$. In this case,  $\mu$ and the integration constant $d$ have a convergent  Taylor series: 
$$
\left\{ \matrix{
  \mu &=& \sum_{n=0}^\infty \mu_n \nu^{2n}
\cr
\cr
d&=& \sum_{n=0}^\infty d_{2n+1}\nu^{2n+1}
}
\right.,~~~\nu\to0
$$
The first terms are 
$$
\mu_0=-{1\over 2},~~~\mu_2=-{\pi\sqrt{3}\over 2},~~~\mu_4=-{\sqrt{3}\pi^3\over 8},~~~\mu_6=-{17\over 240}\pi^5\sqrt{3}
$$ 
and
$$
d_1= {i\pi\over 2}-4\ln(2)-{\pi\sqrt{3}\over 2},~~~d_3=0,~~~
d_5=-\left[{3~\zeta(3)\over 2}+{\pi^3\sqrt{3}\over 30}\right]\pi^2,
$$
$$
d_7=
\left[3~\zeta(5)-{3\pi^2~\zeta(3)\over 4}-{83\pi^5\sqrt{3}\over 7560}\right]\pi^2.
$$
\epr  

\noindent
{\it Proof:} Section \ref{proofchazyd}. 

\vskip 0.2 cm 
The analogous of Proposition \ref{seki} holds:

\bpr
Consider the branch (\ref{pro2y}) associated to the monodromy data $\mu\in{\bf R}$ and $
  p_{1x}=p_{01}=p_{0x}=-2\cosh(2\pi \nu)
$. 
If $\nu\to 0$, then $\mu \to -{1\over 2}$ and,  for $x\neq 0$, there exists: 
$$
\lim_{\nu\to 0} y_1(x)= P_1^{(-1/2)}(\ln x)
$$
where $
 P_1^{(-1/2)}(\ln x)
=-\left(\ln x +d_1+2\right)(\ln x+d_1)$ and $
 d_1= {i\pi\over 2}-4\ln(2)-{\pi\sqrt{3}\over 2}
$. 
\epr

\noindent
{\it Proof:}  $\mu\to -{1\over 2}$ because of Proposition \ref{chazyd}.  Substitute  into (\ref{firstordermu}) the series of $d$ and $\mu$ of Proposition \ref{chazyd} and expand 
$
  \exp\{2id\}=1+2id_1\nu+o(\nu)$ and $x^{2i\nu}=1+2i\nu\ln x+O(\nu^2)
$. 
The structure of the coefficients  $A_{1,-1}$, $A_{10}$, $A_{11}$ allows simplification of the divergences $\nu^{-2}$ and $\nu^{-1}$ contained in the coefficients themselves. Therefore, $y_1(x)$ is expanded in series for $\nu\to 0$, and by direct computation it is easily verified that $
  y_1(x)\to  P_1^{(-1/2)}(\ln x)$ when $\nu \to 0$.

\rightline{\qed}

We verified, up to $n=4$, that every  $y_n(x)$ converges, for $\nu \to 0$, to a polynomial $P_n^{(-1/2)}(\ln x)$. We  conjecture again that if  $\mu=- {1\over 2}$, to the monodromy data $p_{0x}=p_{x1}=p_{01}=-2$ a branch $y(x)$  is associated with asymptotic behavior
\be
\label{trieste}
  {1\over y(x)}\sim \sum_{n=1}^\infty x^{n-1}P_n^{(-1/2)}(\ln x),~~~x\to 0
\ee
where 
$$
P_n^{(-1/2)}(\ln x)= \sum_{N=0}^{2n}p_N~\ln^N(x),~~~~~p_N\in{\bf C}
$$ 
The conjecture is true. It is well known that when $\mu=- {1\over 2}$, to the monodromy data $p_{0x}=p_{x1}=p_{01}=-2$ a one parameter class of  {\it Chazy solutions } of $PVI_{-{1\over 2}}$ is associated. The result is established in \cite{MazzChazy}. Such solutions form a one parameter class,  which includes  (\ref{trieste}) (no parameter in  (\ref{trieste})). 
Therefore, the limit of (\ref{pro2y}) for $\nu \to 0$, $p_{0x}=p_{x1}=p_{01}$ and $\mu $ real, is one element in the class of Chazy solutions. Similar result is estabilished in Section 3.1, Lemma 9, of  \cite{MazzChazy}.

%%%%%%%%%%%%%%%%%%%%%%%%%%%%%%%%%%%%%%%%%%%%%%%%%

\section{Example of Picard Solutions}
\label{exempicard}
Picard solutions \cite{Picard} occurr for $\mu={1\over 2}$. Their example  shows that, when we consider $y(x)$ on the universal covering of the puntured neighborhood of zero, then the poles accumulate at zero along spirals. The Picard solutions of $PVI_{1/2}$ are 
\be
\label{eepi}
  y(x) =\wp\left( \nu_1 \omega_1(x)
  +\nu_2\omega_2(x);\omega_1,\omega^{}_2 \right)+{1+x\over
  3},~~~~~\nu_1,\nu_2\in{\bf C},
\ee
 where the half-periods are  $
\omega_1(x)= {\bf K}(x)$ and $\omega_2(x)=i{\bf K}(1-x)$, and
$$
{\bf K}(x):=\int_0^1 {d\zeta\over \sqrt{
(1-\zeta^2)(1-x\zeta^2)
}
}
={\pi\over 2} F\left({1\over 2},{1\over 2},1;x\right)
$$
 A branch is fixed by the cuts $|\arg x|<\pi$, $|\arg(1-x)|<\pi$.  For $|x|<1$ and $|\arg x|<\pi$, we can write 
$$
\omega_2(x)= -{i\over 2}[F\left({1\over 2},{1\over 2},
1;x\right)\ln(x)+F_1(x)],~~~~|\arg x|<\pi 
$$
where
$$
F_1(x)=
\sum_{n=0}^{\infty}{ \left[\left({1 \over 2}\right)_n\right]^2
  \over (n!)^2 } 2\left[ \psi\left(n+{1\over 2}\right) - \psi(n+1)\right]
x^n,~~~~~|x|<1
$$
 Recal that $\psi(z) = 
{d \over dz}\ln \Gamma(z)$, $ \psi\left({1\over 2}\right) = -\gamma -2 \ln
2$, $ \psi(1)=-\gamma$, and $\psi(a+n)=\psi(a)+\sum_{l=0}^{n-1} {1\over a+l}
$.  
 The behavior of $y(x)$ at $x=0$ follows from the Fourier expansion of $\wp$ at $x=0$. When $x\to 0$ and $\nu_2=2i\nu$, $\nu\in{\bf R}$, this is of type (\ref{pro2y}), with $d={\pi\nu_1\over 2}-\nu\ln 16
$.

\begin{figure}
\epsfxsize=9cm

\epsfysize=9cm

\centerline{\epsffile{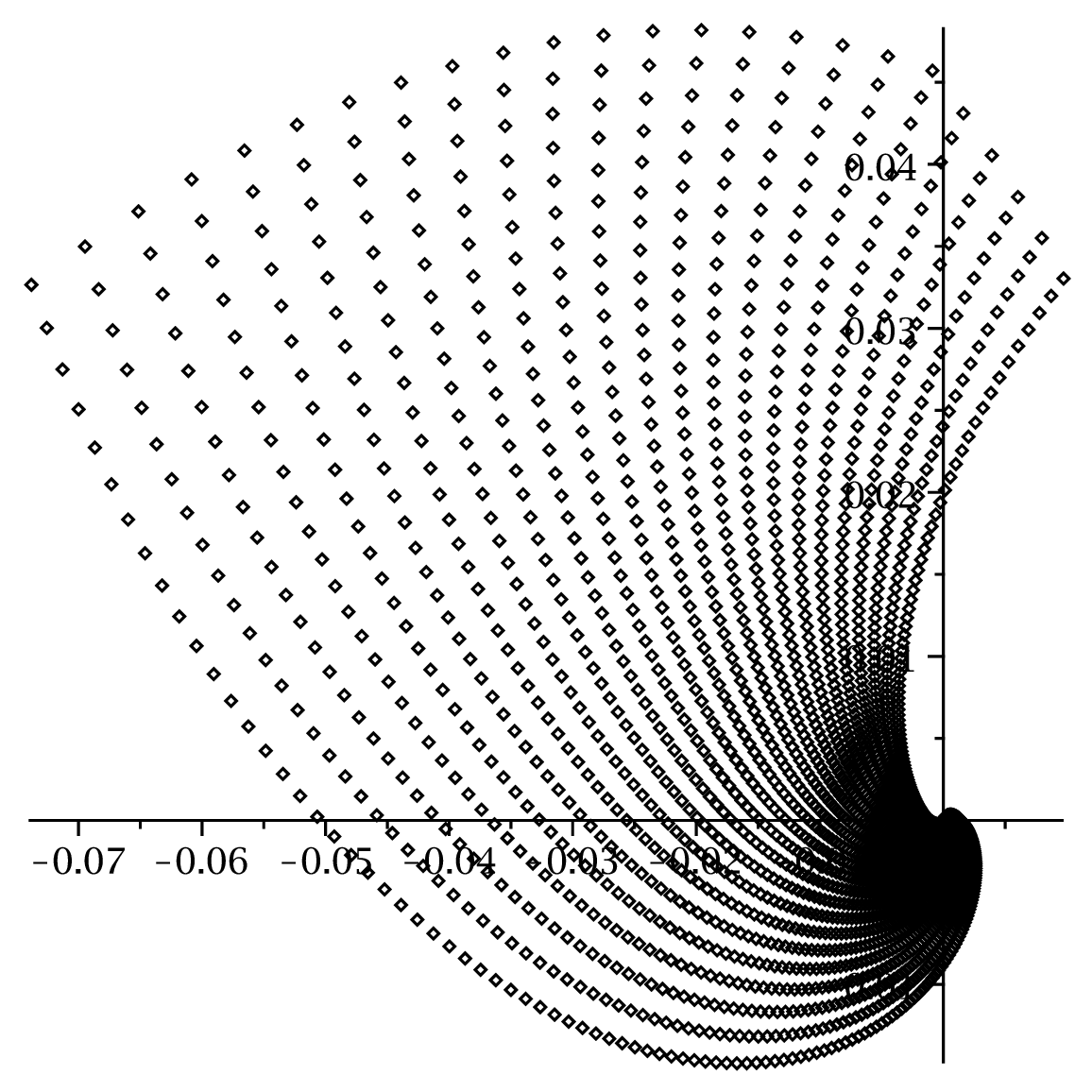}}
\caption{First approximation $x_{kN}$ of the poles of a Picard Solution on the universal covering,  projected onto the $x$ plane.  The figures shows the case  $\nu_1=100-80i$, $\nu_2=1-129i$. Twenty four spirals are displayed, for    $33\leq N\leq 56$. For any $N$, the poles are along a spiral, accumulating at $x=0$  as $k$ increases. In the picture, $k>190$.}
\label{spirale}
\end{figure}

\begin{figure}
\epsfxsize=9cm

\epsfysize=9cm

\centerline{\epsffile{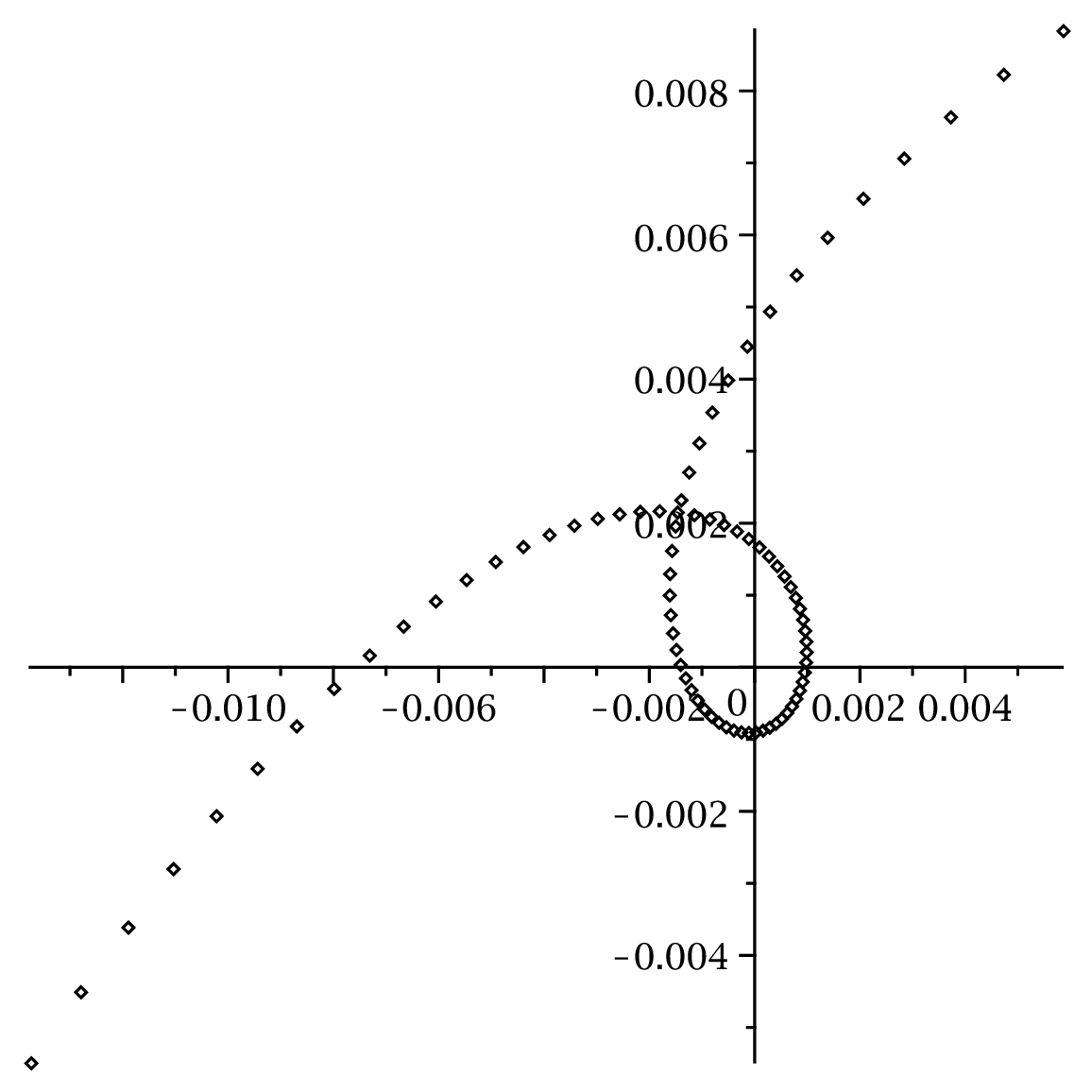}}
\caption{Case of figure \ref{spirale} ($\nu_1=100-80i$, $\nu_2=1-129i$), with fixed $k=250$ and varying $N$ ($-30\leq N\leq 50$).}
\label{plotN}
\end{figure}

 $F_1(x)$ and $F(x)$ are single valued for $|x|<1$, and multi-valuedness of $\omega_2(x)$ comes form $\ln x$.   Thus, $y(x)$ may be regarded as defined on   the universal covering of ${\bf C}\backslash\{0,1,\infty\}$, and in particular, for $|x|<1$,  on the universal covering of a punctured neighborhood of $x=0$. 
The poles on the universal covering of a punctured neighborhood of $x=0$ can be determined. They are a double sequence of points $\xi_{kN}$, solutions of the  equation $
 \nu_1\omega_1(x)+\nu_2\omega_2(x)=
2k\omega_1(x)-2N\omega_2(x)$, $k,N\in{\bf Z}
$, namely: 
$$
 {\nu_2+2N\over 2i} \ln{x\over 16} +{\pi \nu_1\over 2} + {\nu_2+2N\over 2i} 
\left[{F_1(x) \over F(x)}+\ln 16\right]=k\pi,~~~|x|<1,
$$
where ${F_1\over F}+\ln16$  vanishing as $x$, when $x\to 0$. Thus, one can write a pole as 
$$
  \xi_{kN}=x_{kN}\Bigl(1+\delta\left(x_{kN}\right)\Bigr)
$$
 where $x_{kN}$ solves 
$$
 {\nu_2+2N\over 2i} \ln\left({x_{kN}\over 16}\right) +{\pi \nu_1\over 2} =k\pi,~~~k\in{\bf Z}
~~~
\Longrightarrow
 ~~~ x_{kN}= 16 \exp\left\{i\pi{2k-\nu_1\over 2N+\nu_2} \right\}
$$
 The correction  $\delta(x_{kN})$ is expected to vanish if $x_{kN}\to 0$. It is solution of the equation 
\be
\label{secondaeq}
\ln\left( 1+\delta(x_{kN}) \right)+{F_1\Bigl( x_{kN}\left(1+\delta(x_{kN})\right) \Bigr)
\over
 F\Bigl( x_{kN}\left(1+\delta(x_{kN})\right) \Bigr)} +\ln 16 = 0
\ee
One needs to notice the following facts:

\vskip 0.2 cm 
 1) The $x_{kN}$'s can be written as:
 \be
\label{zerosPP}
x_{kN}=x_{0N}\exp\left\{2\pi~{\Im\nu_2\over |\nu_2+2N|^2}~k\right\}\exp\left\{ 2\pi i~ {\Re\nu_2+2N\over   |\nu_2+2N|^2}~k \right\},
\ee
$$
 \hbox{where } ~x_{0N}= 16 \exp\left\{-{i\pi \nu_1\over \nu_2+2N}
\right\}
$$
In order to ensure that $|x_{kN}|<1$, 
the sign of $k$ is chosen $\hbox{sgn}(k)=-\hbox{sgn}(\Im \nu_2)$ and, for given $N$, $|k|$ is sufficiently big.   The above form (\ref{zerosPP}) makes it clear that {\it for any fixed $N$, the $x_{kN}$ lie along a spiral, accumulating at $x=0$ as $k\to \infty$}. The index $N$ singles out the spiral, while $k$ gives the dynamics of the $x_{kN}$'s along that spiral. See figure   \ref{spirale}. 
 On the other hand, if $k$ is fixed and $N$ varies, the typical distribution of the $x_{kN}$'s  is in figure \ref{plotN}. For fixed $k$, $N$ cannot be too big, otherwise $|x_{kN}|$ becomes greater then 1 (and tends to 16 as $N\to\pm \infty$), in which case the local analysis makes no  sense.  
\vskip 0.2 cm 
2) If $\Re \nu_2=0$, no $N$-spiral is a ray. If $\Re \nu_2=2l$, $l\in{\bf Z}$, the spiral for $N=-l$ is  a ray. For example, let $\nu_2=2i\nu$, $\nu\in{\bf R}$. Then, the $x_{k0}$'s are $
 x_{k0}=16\exp\left\{ -{\pi\nu_1\over 2\nu}  -{k\pi\over \nu}
\right\}  
$. 
They  lie along the ray of angle $-{\pi\Im\nu_1\over 2\nu}$, in a disk around $x=0$ of radius less than 1, provided that $k>{\nu\over \pi}\ln16 -{1\over 2} \Re \nu_1$. They are an example of the zeros  of Theorem \ref{theo1}.

\vskip 0.2 cm 
3) Fix the branch cut $-\pi \leq \arg x <\pi$.  Every spiral leaves the cut neighborhood of $x=0$ as $|k|$ increases, and this eventually happens for any $N$. Thus, there are no $x_{kN}$ in a sufficiently small neighborhood (with branch cut) of $x=0$, except possibly in the case  $\Re \nu_2=2l$, $l\in{\bf Z}$, when  the $x_{k,-l}$'s lie along a ray.

\begin{figure}
\epsfxsize=8cm

\epsfysize=8cm

\centerline{\epsffile{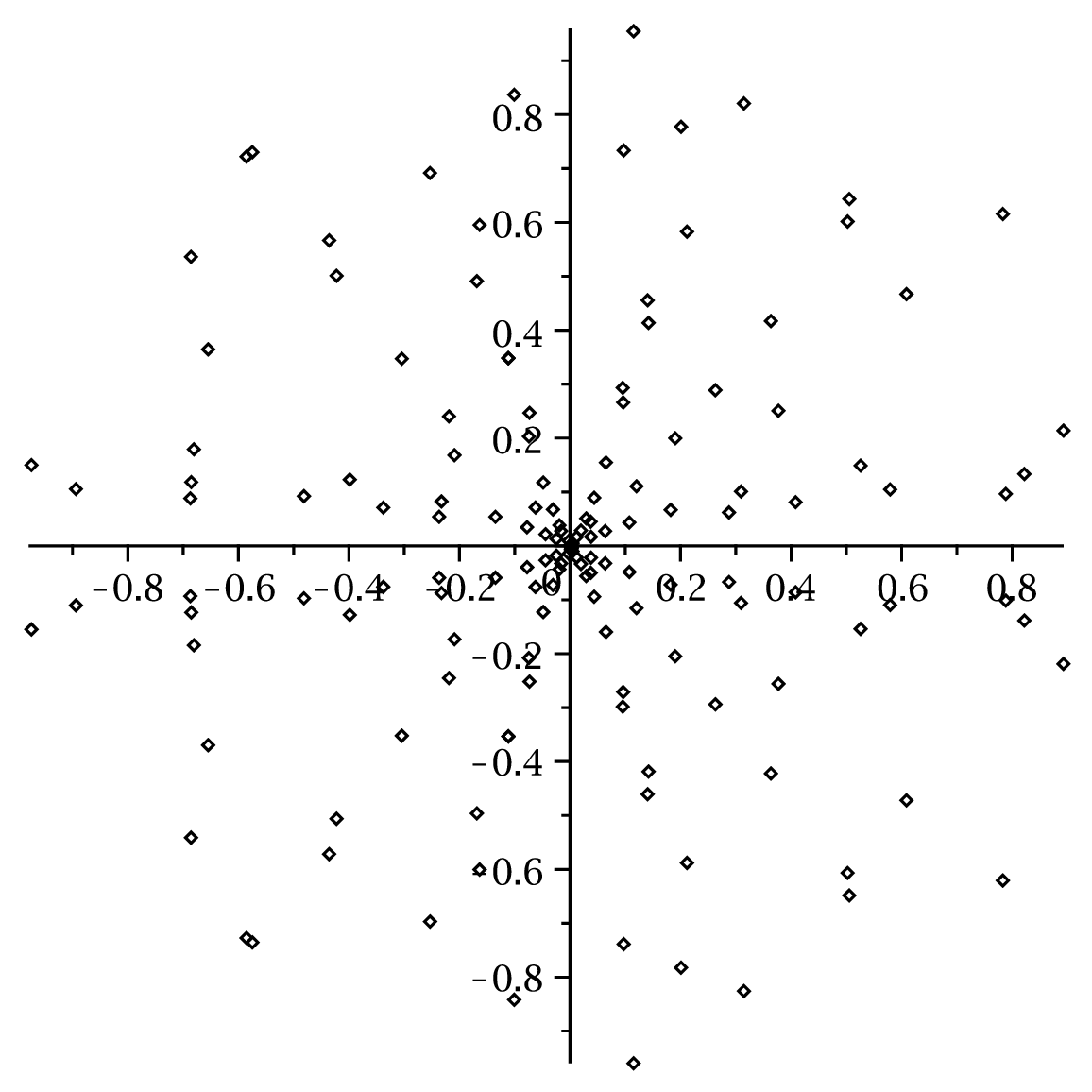}}
\caption{Figure shows the projection on the $x$-plane of the finite number of $x_{kN}$'s, $N\neq 0$,  which lie in $|x|<1$ with $-25\pi<\arg x<23\pi$, for $\nu_1=\nu_2=i/3$. They exist for  $-4\leq N\leq 4$, $N\neq0$ (for $N=0$ there is an infinite sequence $x_{k0}$  on the negative real axis, accumulating at $x=0$, not depicted here).}
\label{poldomain}
\end{figure}

\begin{figure}
\epsfxsize=8cm

\epsfysize=8cm

\centerline{\epsffile{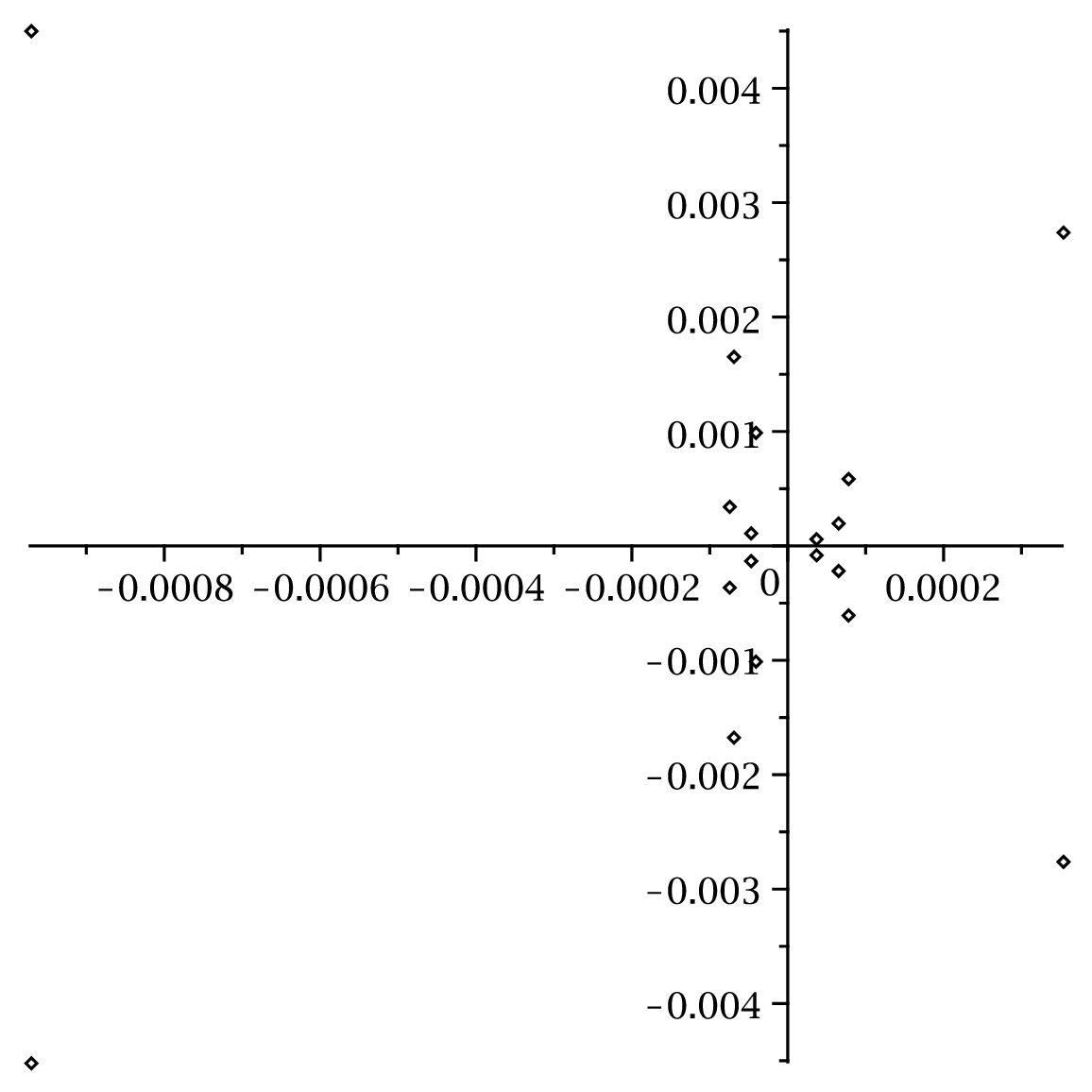}}
\caption{Zoom of figure \ref{poldomain}. At the scale in the figure, only $N=-1$ and $N=1$ appear. $N=0$ is not represented. For $N\neq 0$, there are no  $x_{kN}$ in a neighborhood  sufficiently small of $x=0$. }
\label{poldomain1}
\end{figure}

\vskip 0.3 cm 
\noindent
{\bf Example:} We consider  the case  $\nu_1=\nu_2=i/3$. First, we find the zeros   $x_{kN}$ which lie in a neighborhood of $x=0$, with radius less the 1 and branch cut $-\pi \leq \arg x <\pi$.  The values of $N$ and $k$ are determined imposing  $|x_{kN}|<1$,  $-\pi\leq\arg x_{kN}<\pi$. 
 This gives a system of inequalities 
$$\left\{
\matrix{ 
k < -{1\over 3\pi}(72\ln(2)N^2+2\ln(2)+3\pi N)
\cr
\cr
 k \leq -N, ~~~{18N^2+1\over 18 N} < k
}
\right.~~~~~\hbox{if }N<0
$$
or
$$  \left\{
\matrix{ 
k < -{1\over 3\pi}(72\ln(2)N^2+2\ln(2)+3\pi N)
\cr
\cr
k<{18N^2+1\over 18 N},~~~-N\leq k
}
 \right.~~~~~\hbox{if }N>0
$$
or
$$
k < -{2\ln(2)\over 3\pi}=-0.147...~~~~~\hbox{if }N=0
$$
The above are 
satisfied  only for $N=0$ and $k\leq -1$, which means that only the ray occurring for $N=0$ (the negative real axis) is allowed.  
\begin{figure}
\epsfxsize=9cm

\epsfysize=9cm

\centerline{\epsffile{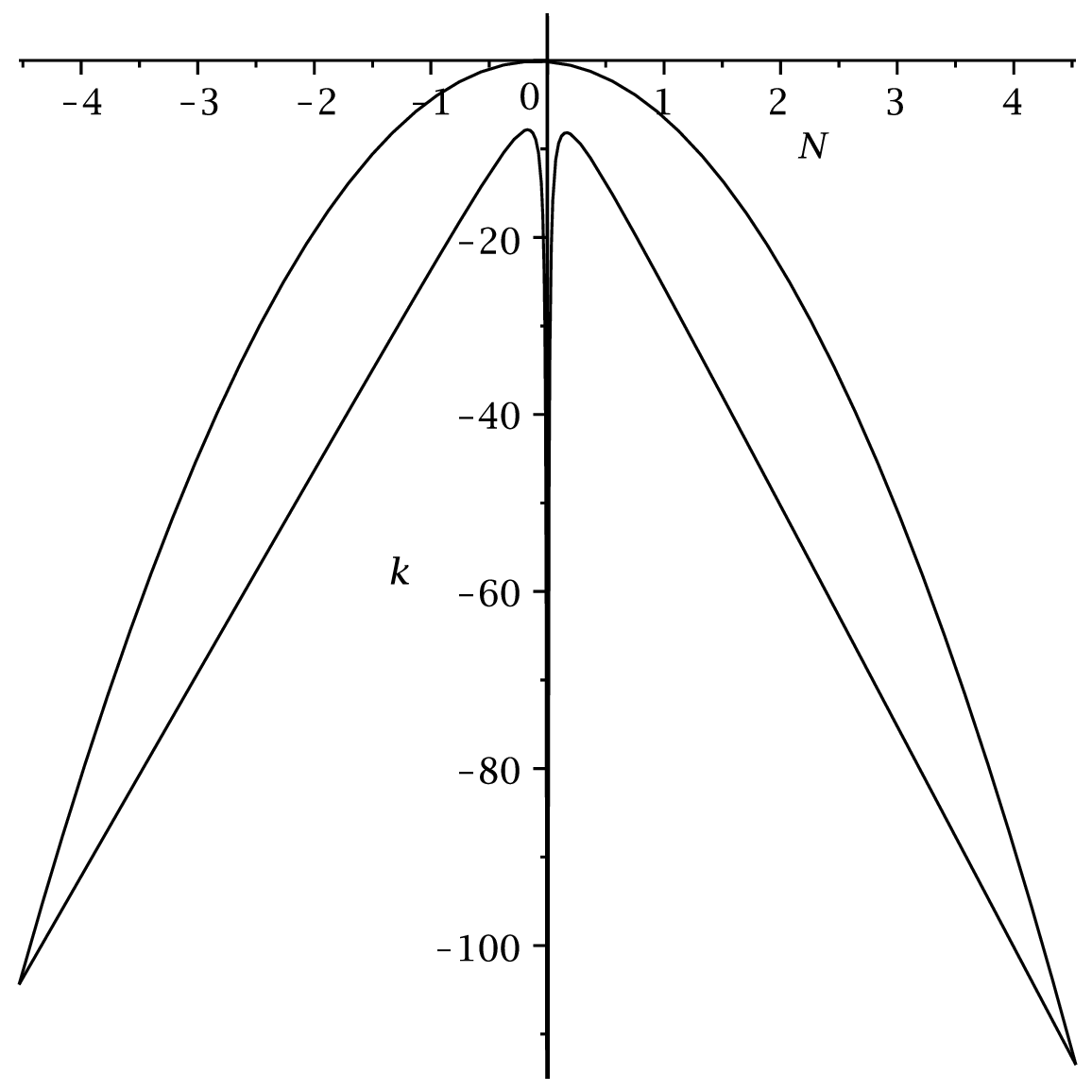}}
\caption{$\nu_1=\nu_2=i/3$. $(N,k)$ plane. When $N$ and $k$ take values in the region delimited by the curves, the points $x_{kN}$ lie in $|x|<1$ with $-25\pi\leq \arg x<23\pi$.  Allowed values of $N$ are $-4\leq N\leq 4$. For $N=0$, any $k\leq-1$ is contained in the region ($N=0$ is a vertical  asymptotic line for the curves).}
\label{domain}
\end{figure}
Then,  we consider a portion of the universal covering, by imposing that $|x_{kN}|<1$ and  $-25\pi<\arg x_{kN} <23\pi$. Again, one obtains inequalities, graphically represented in figure \ref{domain}. The points  $(N,k)$ satisfying the inequalities are   inside the region bounded by the curves of figure \ref{domain}. Therefore, the points $x_{kN}$ which  lie in $\{x~|~|x|<1,~-25\pi\leq\arg x<23 \pi\}$ are: the infinite sequence of  poles on the negative real axis, corresponding to $N=0$ and $k\leq -1$, plus  {\it only a finite number of poles},   corresponding to  $-4\leq N\leq 4$, $N\neq 0$, and $k$ inside the region bounded by the curves of figure \ref{domain}.   They are represented in figures \ref{poldomain} and \ref{poldomain1}.


\vskip 0.3 cm
We now  analyze equation (\ref{secondaeq}). It can be rewritten as
\be
\label{Lagra}
x_{kN}\Bigl(1+\delta(x_{kN})\Bigr)=x_{kN}~{1\over 16} \exp\left\{
-{F_1\Bigl(x_{kN}( 1+\delta(x_{kN}) )\Bigr)\over F
\Bigl(x_{kN}(1+\delta(x_{kN})  )\Bigr)}
\right\}
\ee
Droping the indeces $k$ and $N$ and letting $\xi=x(1+ \delta(x))$, the above  is the equation
\be
\label{Lagra1}
  \xi=x~{1\over 16} \exp\left\{ -{F_1(\xi)\over F(\xi)}\right\}
\ee
 Lagrange inversion theorem can be applied to  (\ref{Lagra1}), because ${1\over 16} \exp\left\{ -F_1(\xi)\over F(\xi)\right\}$ is analytic inside a disk of radius $\xi_0 < 1$, centered at $\xi=0$. The condition on $x$ in order to apply the theorem is $
 |x|<\Bigl|16~\xi\Bigr|\left|\exp\left\{ F_1(\xi)\over F(\xi)\right\}\right|
$,  
when $\xi$ is on the contour of the disk. For such $x$, Lagrange theorem
%%%
%%
%
\footnote{{\it Lagrange Inversion Theorem} (Lagrange 1770 -- see Whittaker \& Watson, a Course of Modern Analysis, pag 133): Let $\phi(z)$ be analytic on and inside a contour $C$ surrounding a point $a$, and let $x$ be such that 
$$
|x~\phi(z)|<|z-a|
$$
at all points  $z\in C$. Then, the equation
$$
x~\phi(\xi)=\xi-a
$$
has one root in  the interior of $C$:
$$
  \xi(x)=a+\sum_{n=1}^\infty {x^n\over n!}{d^{n-1}\over da^{n-1}}\left[\phi(a)^n\right].
$$
}
%
%%
%%%
 says that (\ref{Lagra1}) has a root with the following  convergent series: 
\be
\label{Lagraexp}
\xi(x)= \sum_{n=1}^\infty {x^n\over n!} ~\left.{d^{n-1} \over da^{n-1}}\left[{1\over 16^n} \exp\left\{ -{n~F_1(a)\over F(a)}\right\}\right]~\right|_{a=0}
\ee

$$
=x-{1\over 2}x^2+{11\over 64}x^3-{3\over 64}x^4+{359\over 32768}x^5-{75\over 32768}x^6+{919\over 2097152}x^7+...
$$
The points $
 \xi_{kN}:= \xi(x_{kN})
$, for $x_{kN}$ small, are true poles of $y(x)$. 

\vskip 0.2 cm 

The full description of the poles distribution on the universal covering of ${\bf C}\backslash \{0,1,\infty\}$  for the Hitchin solutions is given in \cite{Br}, where two sequences of poles are determined. Hitchin solutions are solution of the Painlev\'e VI equation with coefficients $\alpha=\beta=
\gamma=1/8$ and $\delta=3/8$, and their image through an Okamoto's transformation (see \cite{Okamoto} and also \cite{MartaBoris}) is the Picard solutions. This transformation annihilates one sequence of poles, and conserves the other, which is given by a simple formula in terms of Theta functions. For $x$ small, this formula coincides with the $\xi_{kN}$ determined here by local analysis.

%%%%%%%%%%%%%%%%%%%%%%%%%%%

\section{The General PVI}
\label{generalPVI}
According to \cite{guz2010}, the general PVI 
 admits solutions with expansion (\ref{pro2y}) and coefficients 
$$
 y_n(x)= \sum_{m=-n}^nA_{nm}(\alpha,\beta,\gamma,\delta) e^{2imd}x^{2im\nu}
$$
The coefficients, as algebraic functions of $\alpha,\beta,\gamma,\delta$, can be computed by the procedure of \cite{guz2010}. For example: 
$$ 
A_{1,\pm 1}=
{\alpha+2\nu^2\mp \sqrt{\gamma~(2\alpha-4\nu^2-\gamma)}
\over 8\nu^2},~~~
A_{10}=-{\alpha-2\nu^2-\gamma\over 4\nu^2}.
$$
(Note that it is allowed the freedom $A_{11}\mapsto c A_{11}$ and   $A_{1,-1}\mapsto c^{-1} A_{1,-1}$, $c\in {\bf C}\backslash \{0\}$, which is equivalent to a redefinition of $d$. But $c e^{2id}$ is fixed by the monodromy data). Theorem \ref{theo1} holds, namely 
 $y_1(x)$ has two infinite sequences of zeros which accumulate at $x=0$ along rays, accoding to the formula 
$$
x_k(j)=\exp\left\{
-{d\over \nu} -{i\over 2\nu} \ln \left[
(-)^j \sqrt{{A_{10}^2\over 4 A_{11}^2}-{A_{1,-1}\over A_{11}}}-{A_{10}\over 2 A_{11}} 
\right]
-{k\pi\over \nu}\right\}
,~~~k\in{\bf N},~~~j=1,2 
$$
The argument of the logarithm is fixed once and for all.  The poles of   $y(x)$  asymptotically approach these zeros, as their absolute value tends to zero, with a series of the form 
$$
\xi_k(j)=x_k(j) +\sum_{N=2}^\infty\Delta_N(j) x_k(j)^N, ~~~k\to +\infty,~~~x_k(j)\to 0.
$$ Also Theorem \ref{theo3} holds.

%%%%%%%%%%%%%%%%%%%%%%%%%%%%%%%%%%%%%%%%%%%%%%%%%%%%%%%%%%%%%%%%%%%%%%%%%%%%%%%%
%%%%%%%%%%%%%%%%%%%%%%%%%% PROOF OF THE RESULTS %%%%%%%%%%%%%%%%%%%%%%%%%%%%%%%%%%%
%%%%%%%%%%%%%%%%%%%%%%%%%%%%%%%%%%%%%%%%%%%%%%%%%%%%%%%%%%%%%%%%%%%%%%%%%%%%%%%%%%%%%%
%%%%%%%%%%%%%%%%%%%%%%%%%%%%%%%%%%%%%%%%%%%%%%%%%%%%%%%%%%%%%%%%%%%%%%%%%%%%%%%%%%%%%%%%%

\section{Proof of Theorem \ref{theo1}}
\label{proofpro7}

The formula for the zeros $x_k(j)$ is proved by solving  
$$
y_1(x)=
          {(2\mu-1+2i\nu)^2e^{-2id}x^{-2i\nu}\over 16 \nu^2}-{(2\mu-1)^2-4\nu^2\over 8\nu^2}
+
{(2\mu-1-2i\nu)^2e^{2id}x^{2i\nu}\over 16\nu^2}
=0
$$

 Let $x_k(j)$ be one of the zeros so obtained. In Theorem \ref{theo3} (whose prove is independent of the following and can be done first)  we prove that there  exists $\epsilon>0$ small such that there are no zeros in ${\cal U}(| x_k(j)|,\epsilon)$ of figure \ref{disk}. Thus, it makes sense to look for the zero  $\xi_k(j)$ of $y(x)$ closest to   $x_k(j)$. Because $\epsilon$ is proportional to $| x_k(j)|$, the lenghts of the arc between $x_k(j)$ and $x_k(j)e^{\pm i\epsilon}$ is proportional to $|x_k(j)|^2$. So, one expects that  $\xi_k(j)= x_k(j)+O(x_k(j)^2)$. Therefore, we look for a zero of the form:
$$
 \xi_k(j)= x_k(j)+\Delta,~~~~~\Delta=\Delta(x):=\sum_{n=1}^\infty \Delta_{n+1}x^{n+1},~~~\Delta_{n+1}\in{\bf C}.
$$ 
 Let $Y(x)= \sum_{n=1}^\infty x^{n-1} y_n(x)$. 
Let $k$ be greater than $k_0$, where $k_0$ is the minimum value such that  $x_{k_0}(j)$ lies in the domain of convergence of (\ref{pro2y}). Impose:
\be
\label{ondina}
 0= Y\Bigl(x_k(j)+\Delta\Bigr)=\sum_{n=1}^\infty \Bigl(x_k(j)+\Delta\Bigr)^{n-1}~
y_n\Bigl(x_k(j)+\Delta\Bigr)
\ee
The series converges if $x_k(j)$ and $x_k(j)+\Delta$ lie in the domain of convergence of (\ref{pro2y}). 
Then, expand:  
$$
  y_n\Bigl(x_k(j)+\Delta\Bigr)= \sum_{m=0}^\infty {1\over m!}~\left.{d^m y_n\over dx^m}\right|_{x_k(j)} \Delta^m
$$
Observe that $[x_k(1)]^{2i\nu}=e^{-2id}$,  $[x_k(2)]^{2i\nu}=\left({2\mu-1+2i\nu\over 2\mu-1-2i\nu}\right)^2e^{-2id}$.  It follows that the derivative of  $
 y_n(x)= \sum_{m=-n}^m A_{nm}(\nu,\mu)e^{2imd} x^{2im\nu}
$ computed at $x_k(1)$ or $x_k(2)$ has the structure: 
$$
 \left.{d^N y_n\over dx^N}\right|_{x=x_k(j)}={Y_{nN}^{(j)}\over [x_k(j)]^N},~~~~~j=1,2
$$
where $ Y_{nN}^{(j)} $ are constants which {\it do not} depend on $k$. They depend only on $\nu,\mu$. In particular, $Y_{n0}^{(j)}= y_n(x_k(j))$, and $Y_{10}^{(j)}=y_1(x_k(j))=0$.  We have therefore to determine $\Delta=\Delta(x)$ which solves:
$$
 0= \sum_{n=1}^\infty \left[
\left(
x+\Delta
\right)^{n-1}
~
\sum_{m=0}^\infty {Y_{nm}
  \Delta^m
\over 
m! ~x^m}
\right]
= \sum_{n=1}^\infty x^{n-1}\left[
\left(
1+\tilde{\Delta}
\right)^{n-1}
~
\sum_{m=0}^\infty {Y_{nm}\over m!}~
  \tilde{\Delta}^m
\right],~~~~~\tilde{\Delta}={\Delta\over x}
$$
This is similar to a problem of reversion of a series, though it is not in the form which allows to apply Lagrange inversion theorem to find $\tilde{\Delta}(x)$. Nevertheless, the coefficients are computable by putting equal to zero the coefficients of the powers of $x$ in   
the series expansion (omitting  $k$ and $j$):
$$
 0
 = \sum_{n=1}^\infty \left[x^{n-1}
\left(
1+\sum_{l=1}^\infty\Delta_{l+1}x^{l}
\right)^{n-1}
~\left(
\sum_{m=0}^\infty {Y_{nm}\over m!}
 \left( 
\sum_{l=1}^\infty \Delta_{l+1}x^{l}
 \right)^m
\right) 
\right]
$$
$$
=\sum_{n=1}^\infty \left[
Y_{n0}+\sum_{p=1}^n \pmatrix{n-1 \cr p-1} \sum_{m=1}^\infty {Y_{nm}\over m!}\sum_{l_1,...,l_{m+p-1}=1}^\infty \Delta_{l_1+1}...\Delta_{l_{m+p-1}+1} x^{l_1+...l_{m+p-1}}
\right]x^{n-1}
$$
The above series  determines all the $\Delta_n$'s recursively, provided that $Y_{10}=0$, namely provided that $x_k(j)$ is a zero of $y_1(x)$.  The first terms are: 
$$
\Delta_2= -{Y_{20}\over Y_{11}},~~~~~~~
\Delta_3= -{Y_{30}\over Y_{11}}+{2Y_{20}(Y_{2 0}+Y_{2 1})Y_{1 1}-Y_{1 2}Y_{2 0}^2 \over 2Y_{1 1}^3}
$$
 $$
\Delta_4= -{Y_{40}\over Y_{11}}+{1\over 2 Y_{11}^5} \Bigl\{[(6Y_{3 0}+2Y_{3 1})Y_{2 0}+2Y_{2 1}Y_{3 0}]Y_{1 1}^3+
$$
$$  -
2\Bigl(Y_{2 0}^2+(
3Y_{2 1}+{1\over 2}Y_{2 2})Y_{2 0}+Y_{12}Y_{3 0}+Y_{2 1}^2\Bigr)Y_{2 0}Y_{1 1}^2+
$$
$$
+\Bigl((3Y_{1 2}+{1\over 3}Y_{1 3})Y_{2 0}^3
+3Y_{1 2}Y_{2 0}^2Y_{2 1}\Bigr)Y_{1 1}-Y_{1 2}^2Y_{2 0}^3
\Bigr\}
$$
The above formulas do not hold for double zeros, namely $Y_{11}=0$,  which occur for $\mu={1\over 2}$. In this case, the procedure of inversion of the series works if also $Y_{20}=0$, which is true when $\mu={1\over 2}$. The formulas obtained are more complicated and will be omitted here. The case $\mu={1\over 2}$ corresponds to a sub-class of Picard solutions, and it is solved in Section \ref{exempicard}. 

\vskip 0.2 cm 
\noindent
{\bf Note:} 
The computation, rather complicated, can be done on a computer. To determine $\Delta_2$, $\Delta_3$,... 
$\Delta_{N+1}$ one can write the following expression for a program of symbolic computation like Maple or Mathematica: 
$$
 0= \sum_{n=1}^{N+1} \left[
\left(
x+\sum_{l=1}^{N}\Delta_{l+1}x^{l+1}
\right)^{n-1}
~\left(
\sum_{m=0}^N {Y_{nm}
 \left( 
\sum_{l=1}^{N} \Delta_{l+1}x^{j+1}
 \right)^m
\over 
m! ~x^m}
\right) 
\right]
$$
It is enough to compute  $y_1$ up  the $N$-th derivatives, $y_2$ up to order $(N-1)$-th derivative,...,  $y_{N-1}$ up to first derivative, and $y_N$.

\vskip 0.2 cm 

\noindent
 We  compute the $Y_{nN}^{(j)}$'s from the $A_{nm}$'s, and find  $
 Y_{11}^{(1)}=  -Y_{11}^{(2)}= 2\mu-1$,  $Y_{20}^{(1)}=-Y_{20}^{(1)}=\mu-{1\over 2}$. It follows  that  $
 \Delta_2(1)=\Delta_2(2)=-{1\over 2}$. 
Higher order computations, up to $N=3$, give the other formulas of $\Delta_N(j)$ stated in  the Theorem.   
Note that, though the derivation of the result requires $\mu\neq {1\over 2}$, these formulas  have limit for $\mu\to {1\over 2}$, which coincide with the result of the example of Section \ref{exempicard}.

\vskip 0.2 cm 
 It is necessary to show that the result is consistent. Observe that $\xi_k(j)=x_k(j)-{1\over 2} x_k(j)^2 +O(x_k(j)^3)$. Thus,  we need to show that there exist a $K$ such that for any $k\geq K$   a disk of radius $|x_k(j)|^2$ with center  $x_k(j)$ does not intersect a similar disk around another zero. This would  imply that  for  any $j=1,2$ and for any $k\geq K$, equation  (\ref{ondina}) 
has a unique formal solution $x_k(j)+\Delta^{(j)}(x_k(j))$,  with $
\Delta^{(j)}(x)= \sum_{N\geq2}\Delta_{N}(j)~x^N$,  
where the $\Delta_N(j)$'s have been uniquely constructed by the procedure above. We distinguish two orderings  of the zeros

{\bf --} When the ordering is $|x_{k+1}(1)|<|x_k(2)|<|x_k(1)|$. Two cases must be considered:  i) $x_k(1)$ and $x_{k}(2)$; ii)   $x_k(2)$ and $x_{k+1}(1)$. 
In case i), one must check if the following holds:
$$ 
 |x_{k}(2)|+|x_{k}(2)|^2< |x_{k}(1)|-|x_{k}(1)|^2
$$
Since $|x_{k}(2)|=e^{{\theta\over \nu}} |x_{k}(1)|$, where $-\pi<\theta= \arg\left({2\mu-1+2i\nu\over 2\mu-1-2i\nu}\right)\leq 0$, the above becomes:
$$
|x_{k}(1)|e^{{\theta\over \nu}}+|x_{k}(1)|^2e^{{2\theta\over \nu}}< |x_{k}(1)|-|x_{k}(1)|^2
$$ 
which holds for 
$$
0<|x_k(1)|< {1-e^{{\theta\over \nu}}\over 1+e^{{2\theta\over \nu}}}<1
$$
In case ii), one must check if the following holds:
$$ 
 |x_{k+1}(1)|+|x_{k+1}(1)|^2< |x_{k}(2)|-|x_{k}(2)|^2
$$
Since $|x_{k+1}(1)|=e^{-{\pi+\theta\over \nu}} |x_{k}(2)|$, the above becomes:
$$
|x_{k}(2)|e^{-{\pi+\theta\over \nu}}+|x_{k}(2)|^2e^{-{2(\pi+\theta)\over \nu}}< |x_{k}(2)|-|x_{k}(2)|^2
$$
which holds for 
$$
0<|x_k(2)|< {1-e^{-{\pi+\theta\over \nu}}\over 1+e^{-{2(\pi+\theta)\over \nu}}}<1
$$ 
Therefore, there exist $K$ big enough such that  for any $k\geq K$ the above inequalities i) and ii) are satisfied. 

\vskip 0.2 cm

{\bf --}   When the ordering  is $|x_{k}(1)|<|x_k(2)|<|x_{k-1}(1)|$. Two cases must be considered:  i) $x_{k-1}(1)$ and $x_{k}(2)$; ii)   $x_k(2)$ and $x_{k}(1)$. We proceed and conclude in a similar way that there exist $K$ big enough such that  for any $k\geq K$ the  inequalities are satisfied. 

\rightline{\qed}

%%%%%%%%%%%%%%%%%%%%%%%%%%%%%%%%%%%%%%%%%%%%%%%%%%%%%%%%%

\section{Proof of Theorem \ref{theo3}}
\label{proofpro6}

The question to be answered is where the zeros of  ${1\over y(x)} =\sum_{n=1}^\infty x^{n-1} y_n(x)$ do {\it not} lie. First, note that $x=0$ is an essential singularity for $y_1(x)$ (with $|y_1(x)|$ bounded), and it is not a zero. Let us write: 
$$
  \sum_{n=1}^\infty x^{n-1} y_n(x)= y_1(x)+xf(x)
$$
where $f(x)$ is a bounded function in a disk of radius smaller than the radius of convergence of the series above. Let again $k\geq k_0$, where $k_0$ is such that $x_{k_0}(j)$ is the biggest zero in the domain of convergence. $j$ is either $j=1$ or $j=2$, depending on the order  $|x_{k}(2)|<|x_k(1)|$ or  $|x_{k}(1)|<|x_k(2)|$, the first case occurring when  $-\pi < \arg\left({2\mu-1+2i\nu\over 2\mu-1-2i\nu}\right)\leq 0 $, and the second when  $0 < \arg\left({2\mu-1+2i\nu\over 2\mu-1-2i\nu}\right)\leq \pi $. Boundedness of $f(x)$ means that there exists $C_f>0$ such that  
$|f(x)|\leq C_f$ for any $x$ in the disk $|x|\leq |x_{k_0}(j)|$. 
 
 \vskip 0.2 cm 
First, we prove the statement for the case $\mu=-1$. 
 The zeros of $y_1(x)$ are on the negative immaginary axis and $x=0$. For any $0<\epsilon<\pi$, $y_1(x)$ is
 bounded and not vanishing in the domain ${\cal U}(|x_{k_0}(1)|,\epsilon)$ of the type represented 
in figure  \ref{disk}.  Namely, for any $\epsilon>0$ small, there exist $C_\epsilon>0$ such that   $|y_1(x)|\geq C_\epsilon$ for $x\in {\cal U}(|x_{k_0}(1)|,\epsilon)$. Observe that $C_\epsilon\to 0$ as $\epsilon\to 0$, because $x$ may take values closer  to the zeros of $y_1(x)$ as $\epsilon$ gets smaller. The following simple estimate holds for 
$x\in  {\cal U}(|x_{k_0}(1)|,\epsilon)$: 
\be
\label{ineq}
 \Bigl|y_1(x)+xf(x)\Bigr|\geq \Bigl||y_1(x)|-|x||f(x)| \Bigr|\geq |y_1(x)|-|x|C_f \geq C_\epsilon-|x|C_f
\ee
We have:
$$
C_\epsilon-|x|C_f >0~~~\Longleftrightarrow~~~ |x|<{C_\epsilon\over C_f}
$$
 Namely, $y_1(x)+xf(x)$ has no zeros in  ${\cal U}(R_\epsilon,\epsilon)$, for any $0<R_\epsilon <{C_\epsilon\over C_f}$.

\vskip 0.2 cm
In order to prove the  estimate $
C_\epsilon= 3 \tan(\epsilon)~ (1+O(\tan (\epsilon)^2))$ for $\epsilon\to 0$, 
\begin{figure}
\epsfxsize=9cm

\epsfysize=9cm

\centerline{\epsffile{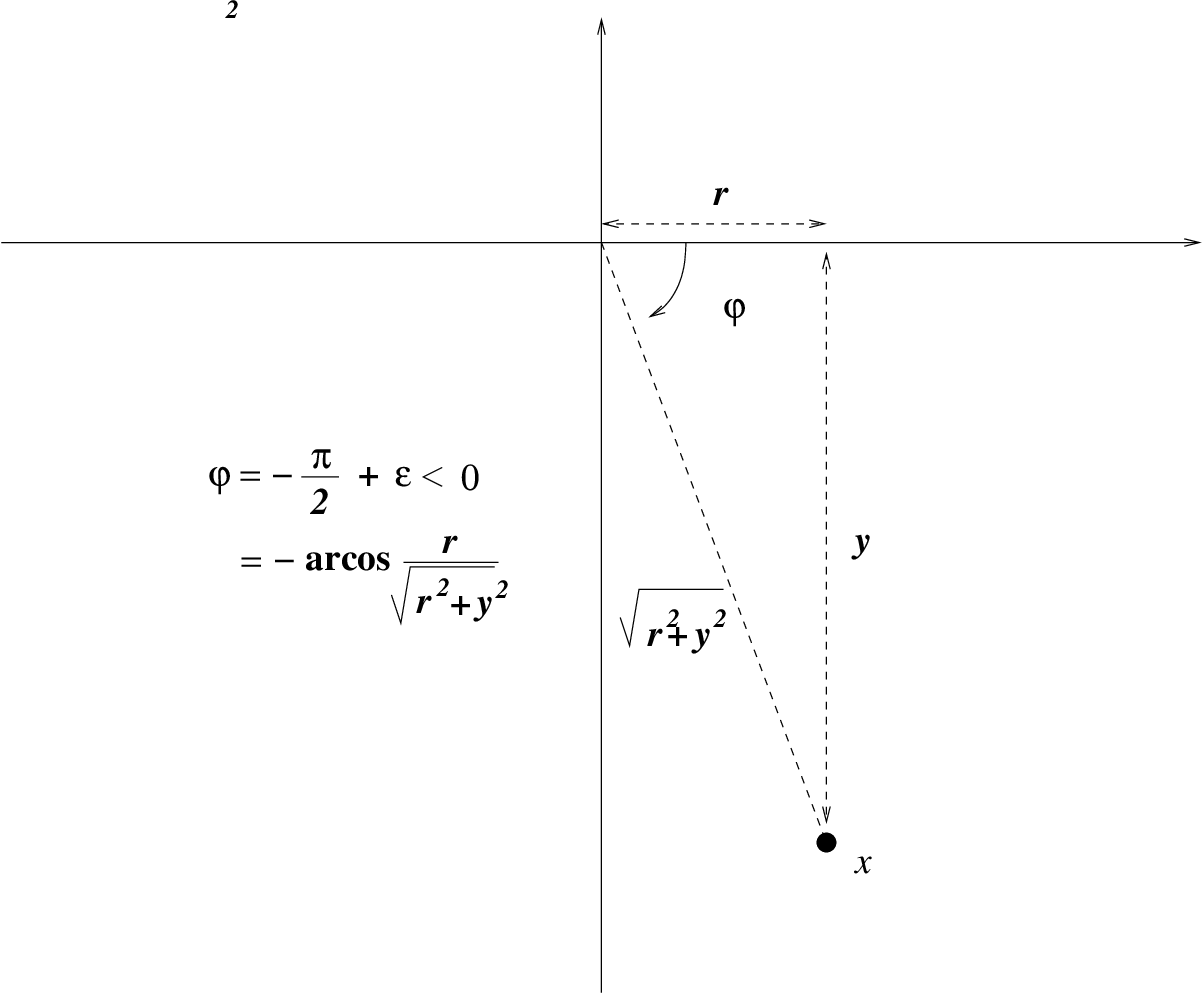}}
\caption{}
\label{rectangle}
\end{figure}
  we evaluate $|y_1(x)|$ along the two rays of angle $-{\pi\over 2}\pm \epsilon$, which are the boundary of the sectorial cut of $ {\cal U}(|x_{k_0}(1)|,\epsilon)$.  On these rays  $|y_1(x)|$ take minimal values as $\epsilon \to 0$, because the zeros on the negative immaginary axis are approached. We write
$$ 
y_1(x)= {(3-2i\nu)^2\over 16\nu^2}e^{-2id}x^{-2i\nu} -{9-4\nu^2\over 8\nu^2} +{(3+2i\nu)^2\over 16\nu^2} e^{2id}x^{2i\nu}
$$
We consider the ray of angle $-{\pi\over 2}+ \epsilon$, the other being analogous.
Let us write $(3\pm2i\nu)^2= (9+4\nu^2)\exp\{\pm 2i\psi\}$, where 
$\psi=\hbox{arcos}\left({3\over\sqrt{9+4\nu^2}} \right)>0$.  Let us also write $x=r-iw$, $w>0$ along the ray (see figure \ref{rectangle}). 
Thus: 
$$ (r-iw)^{2i\nu} 
= e^{2\nu \arccos\left({r\over \sqrt{r^2+w^2}}\right)} ~e^{i\nu \ln(r^2+w^2)},
$$
where $\arccos\left({r\over \sqrt{r^2+w^2}}\right)={\pi\over 2}-\epsilon>0$. Recall that $\Im d = {\pi \over 2} \nu$. It follows that:  
\be
\label{star}
|y_1(x)|^2 = \left({9+4\nu^2\over 8\nu^2}\right)^2 \left\{
\left(\cosh(a)\cos(b)-{9-4\nu^2\over9+4\nu^2}\right)^2+ \sinh(a)^2\sin(b)^2
\right\}.
\ee
$$
a= \nu\left(\pi -2 \arccos {r\over \sqrt{r^2+w^2}}\right),~~~~~b= \nu\ln(r^2+w^2)+2(\psi+\Re d).
$$
When $\epsilon\to 0$, then $r\to 0$. The minimal values are obtained when $w= |x_k(j)|$,  $k\in {\bf N}$, $j=1,2$, where $x_k(j)$  is a zero of $y_1(x)$. To this end, before proceeding, observe that if we let $r\to 0$ with $w\neq 0$ fixed, then $-\arccos{r\over \sqrt{r^2+w^2}}\to -{\pi\over 2}$ and $a\to 0$. Therefore:
$$
 |y_1(x)|^2 \to \left({9+4\nu^2\over 8\nu^2}\right)^2 \left(\cos b - {9-4\nu^2\over 9+4\nu^2} \right)^2,~~~
 b=2\nu \ln w +2(\psi +\Re d).
$$
This takes minimal value $=0$ if:
$$
 b=\arccos{9-4\nu^2\over 9+4\nu^2}-2k\pi,~~~~~k\in {\bf N},
$$
namely: 
 \be
\label{maldispalla}
 w= \exp\left\{-{\Re d\over \nu} -{\psi\over \nu} +{1\over 2\nu} \arccos {9-4\nu^2\over 9+4\nu^2}\right\}~\exp\left\{-k{\pi\over \nu}\right\}.
\ee
We prove that the above coincides with the zeros of $y_1(x)$. To see it, observe that:
$$
 0<\psi=\arccos{3\over \sqrt{9+4\nu^2}}=-i \ln\left({3+2i\nu\over \sqrt{9+4\nu^2}}\right)
$$
$$
\arccos {9-4\nu^2\over 9+4\nu^2} = -i \ln\left({9-4\nu^2\over 9+4\nu^2}\pm \sqrt{\left({9-4\nu^2\over 9+4\nu^2}\right)^2-1  }\right)
= -2i\ln\left( {3\pm 2i\nu \over \sqrt{9+4\nu^2}} \right)
$$  
which implies: 
$$
-{\psi\over \nu} +{1\over 2\nu} \arccos {9-4\nu^2\over 9+4\nu^2}= -{i\over \nu} 
\ln\left( {3\pm 2i\nu\over 3+2i\nu} \right)
$$
and then (\ref{maldispalla}) corresponds to  the (absolute value of the) zeros: 
$$
 |x_k(j)|=\exp\left\{-{\Re d \over \nu} + {j-1\over \nu} \arg{3-2i\nu\over 3+2i\nu}\right\}~
\exp\left\{-k{\pi\over \nu}\right\},~~~~~j=1,2
$$
Now we proceed to the estimate of $y_1(x)$ along the radius of angle $-{\pi\over 2}+\epsilon$, when $r\to 0$. We have: 
$$
 0<\arccos{r\over \sqrt{r^2+w^2}}= -i \ln\left({r+iw\over \sqrt{r^2+w^2}} \right)= {\pi \over 2} -{r\over w}+O\left({r^3\over w^3}
\right)
$$
Thus $
 a= {2\nu r\over w} + O\left({r^3\over w^3}
\right)
$, 
and 
$$
\cosh(a)= 1+{2\nu^2r^2\over w^2}+O\left({r^3\over w^3}\right),~~~~~\sinh(a)={2\nu r\over w}+O\left({r^3\over w^3}\right)
$$
Moreover:
$$
b= \nu\ln(r^2+w^2)+2(\psi+\Re d)= \nu\left[ \ln(w^2)+\ln\left(1+{r^2\over w^2}\right)\right]+2(\psi+\Re d)
$$
$$
=2\nu\ln w +2(\psi+\Re d)+{\nu r^2\over w^2}+O\left({r^4\over w^4}\right) 
$$
The evaluation will be done when, for $\epsilon$ small and $r$ small, the position of $r-iw$ is close to a zero of $y_1(x)$, namely to a minimum of $|y_1(x)|$. According to the above computations for $r=0$ and $w$ fixed, we have that $r-iw$ is close to a zero of $y_1(x)$ when  $w= |x_{k}(j)|$. Namely: 
$$ 
 2\nu\ln w +2(\psi+\Re d)= \arccos {9-4\nu^2\over 9+4\nu^2}-2k\pi
$$
($j=1,2$ is given by the to signs of arccos).  Therefore, for small $r$ and $w= |x_{k}(j)|$:
$$
\cos(b)= \cos\left(  \arccos{9-4\nu^2\over 9+4\nu^2}+{\nu r^2\over w^2}+O\left({r^4\over w^4}\right) \right).
$$
Now, use the fact that $\cos(\theta+\delta)=\cos\theta-\delta\sin\theta-{\delta^2\over 2}\cos\theta +O(\delta^3)$ when $\delta\to 0$, to obtain: 
$$
 \cos(b)= {9-4\nu^2\over 9+4\nu^2}-{\nu r^2\over w^2} \sqrt{1-\left( {9-4\nu^2\over 9+4\nu^2}\right)^2}+O\left({r^4\over w^4}\right)= {9-4\nu^2\over 9+4\nu^2} -{12\nu^2\over 9+4\nu^2}~ {r^2\over w^2} +O\left({r^4\over w^4}\right) 
$$
$$
0<\sin(b)= \sqrt{1-\cos(b)^2}= {12\nu\over 9+4\nu^2} \left(1-{9-4\nu^2\over 12}~{r^2\over w^2}+O\left({r^4\over w^4}\right)\right)
$$
Composing all the above expansions into (\ref{star}), we finally obtain the minimal values:
$$
 |y_1(x)|= {3r\over w}\left(1+O\left({r^2\over w^2}\right)\right)=3 \tan(\epsilon)~ (1+O(\tan (\epsilon)^2))
$$

\vskip 0.3 cm
We turn to the general case. 
The domain ${\cal U}(|x_{k_0}(j)|,\epsilon)$ is that of figure \ref{disk}. There is no conceptual change in proving the inequality (\ref{ineq}).  The estimate of $C_{\epsilon}$ is done as for the case $\mu=-1$, for both the sector cuts of  ${\cal U}(|x_{k_0}(j)|,\epsilon)$. Technically, it is more complicated, but methodologically the same. The final result is that the minimal values of $|y_1(x)|$ are 
$$
 |y_1(x)|= {|2\mu-1| r \over w} \left(1+O\left({r^2\over w^2}\right)\right)= |2\mu-1|\tan(\epsilon)~ (1+O(\tan (\epsilon)^2))
$$
Thus $
 C_\epsilon= |2\mu-1|\tan(\epsilon)~ (1+O(\tan (\epsilon)^2)),$ $\epsilon\to 0. 
$ 
In the case $\mu={1\over 2}$, the above procedure should be applied expanding $|y_1|$ at least to order $r^2/w^2$, because the zeros of $y_1(x)$ are double ($x_{k+1}(1)=x_k(2)$), and ${d^2 y_1(x_k(j))\over dx^2}=0$. The final result is
$$
 |y_1(x)|= {2\nu^2 r^2 \over w^2} \left(1+O\left({r^3\over w^3}\right)\right)= 2\nu^2\tan(\epsilon)^2~ (1+O(\tan (\epsilon)^3))
$$
\rightline{\qed}

% % %% % % % % % % % %% % % % % % % % % % % %%  % % % % % % %

\section{Proof of Proposition \ref{lem1}:} 
\label{monodromydata}

The parametrization of a branch (\ref{pro2y}) in terms of monodromy data is computed in section 6 of our \cite{guz2010} (see also Proposition 6 of \cite{guz2010}), where we find the notation $\sigma_0=1+2i\nu$ and the result $2\cos\pi\sigma_0=p_{0x}$. Thus, $\nu$ is determined by   $
\cosh(2\pi\nu)= -{p_{0x}\over 2}
$, 
which has solutions: 
$$
 \exp\{\pm2\pi\nu\}= \pm{1\over 2}\sqrt{p_{0x}^2-4}-{p_{0x}\over 2},~~~~~\nu>0
$$
This proves  (\ref{nugen}).  In \cite{guz2010} $d$ is computed, and the result is:  
$$
d= {i\over 2}~\ln\left\{{8i\nu~r\over\Bigl(2\nu+i(1-2\mu)\Bigr)^2}\right\},
$$
where
$$r=
-{16^{2i\nu}~\Gamma({3\over 2}-{\theta_\infty\over 2}-i\nu)^2 ~\Gamma({\theta_\infty\over 2}+{1\over 2}-i\nu)^2 \over 2i~\nu^3~\sinh(2\pi\nu)^2~\Gamma(-i\nu)^4}~
\times~~~~~~~~~~~~~~~~~~~~~~~~~~~~~~~~~~~~~~~~~~~~
$$
$$~~~~~~~~~~~~~~~\times~\left[{1\over 2}
\Bigl(e^{2\pi\nu}p_{x1}-p_{01}\Bigr)\sinh(2\pi\nu)+\Bigl(\cos(\pi\theta_\infty)+1\Bigr)\Bigl(e^{2\pi\nu}+1\Bigr)\right].
$$
This proves (\ref{generald}). 
In the special cases $2\mu=\pm 2i\nu+2m+1$, we have $p_\infty=p_{0x}$. The  affine cubic becomes a polynomial of degree 2 in $p_{01}$, 
with to solutions $p_{01}= 2-2e^{\pm2\pi\nu}+p_{x1}e^{\pm2\pi\nu}$. 
The formulae for $d$ when  $2\mu=2i\nu+1-2m$, $m=-1,-2,-3,...$, and $2\mu=-2i\nu+2m-1$, $m=0,-1,-2,-3,...$,  can be obtained by substitution into the formula of the generic case. This is not possible in other cases, where we need the formulas computed in our \cite{D4}. In \cite{D4} we find the Jimbo's  solution 
$$
 y_{\hbox{\rm\cite{D4}}}(x)\sim ax^{1-\sigma}, ~~~0<\Re\sigma<1
$$
whose coefficient $a$ is computed in terms of monodromy data also  in the special cases $\sigma  \pm 2\mu=2m$. 
Solution (\ref{pro2y}) can be ritten as:  
$$
 y(x)= {1\over A_{1,-1}e^{-2id} x^{1-\sigma}+A_{10}+A_{11}e^{2id}x^{\sigma-1}+O(x)},~~~~\sigma=1+2i\nu
$$
Now, let $\sigma=1+2i\nu-\varepsilon$, $0<\varepsilon\to 0$, and rewrite
$$
y(x)= {A_{11}^{-1}e^{-2id} x^{1-\sigma}\over
1+\hbox{\small higer order corrections } O(x^\varepsilon)}
$$
Then, 
identify $
  A_{11}^{-1}~e^{-2id}= a$   ($a$ is given in theorem 2, page 301, of \cite{D4}), extract $d$ and let $\varepsilon \to 0$. 
 This completes the proof.  (Note that in \cite{D4}, the  notations of \cite{DM} are used, namely  
 $x_0^2=2-p_{0x}$, $x_1^2=2-p_{x1}$ and  $x_\infty^2= 2-p_{01}$). 

\rightline{\qed}

%%%%%%%%%%%%%%%%%%%%%%%%%%%%%%%%%%%%%%%%%%%%%%%%%%%%

% % % % % % % % % % % % % % % %% % % %  % % % %%  %% %  % % % 

\section{Proof of theorem \ref{theo2} and its Corollary}

\subsection{The case $p_{1x}=p_{01}=p_{0x}$}
\label{forCHAZY}

The behavior at $x=0$  of a branch of a $PVI_\mu$-transcendent  is explicitly parameterized by the monodromy data 
$\theta_\infty=2\mu$ (i.e. $p_\infty$), $p_{0x},p_{01},p_{x1}$ to which it is in one-to-one correspondence. 
Let $\sigma$ be defined by $
 2\cos\pi\sigma=p_{0x}$ and  $0\leq \Re \sigma \leq 1
$. Its value determines the critical behaviors as follows:  
\vskip 0.2 cm 
\noindent
-- When  $0<\Re \sigma <1$, namely $p_{0x}\not\in(-\infty,-2]\cup[2,+\infty)  $, the behavior is (see Jimbo \cite{Jimbo}):
$$
 y(x)= ax^{1-\sigma}(1+O(x^\sigma,x^{1-\sigma})),~~~~~~~a=a(\sigma,\theta_\infty, p_{01},p_{x1})\in{\bf C}
$$ 

\noindent
-- When  $\sigma=2i\nu$, $\nu>0$, namely $p_{0x}>2$,  the behavior is (see \cite{Jimbo} and \cite{guz2010}): 
$$
 y(x)=x\left[{(2\mu-1+2i\nu)^2e^{-2id}x^{-2i\nu}\over 16 \nu^2}-{(2\mu-1)^2-4\nu^2\over 8\nu^2}
+
{(2\mu-1-2i\nu)^2e^{2id}x^{2i\nu}\over 16\nu^2}\right]+O(x^2)
$$
where $d=d(\nu,\theta_\infty, p_{01},p_{x1})\in{\bf C}$. 

\vskip 0.2 cm 
\noindent
-- When  $\sigma=1+2i\nu$, namely $p_{0x}<-2$, the behavior is  (\ref{pro2y}).  

\vskip 0.2 cm
 In the special case $
p_{0x}=p_{x1}=p_{01}$, the cubic surface (\ref{cubic}) becomes the curve (\ref{cubic1}), depicted in  figure \ref{CUBIC}  for $\cos\pi\theta_\infty\in{\bf R}$. It has three branches when  $-1<\cos\pi\theta_\infty<1$, namely when $\theta_\infty$ is real. It has  double points for $\cos\pi\theta_\infty=-1$, namely $\theta_\infty= 2m+1$, $m\in {\bf Z}$, and for $\cos\pi\theta_\infty=1$, namely $\theta_\infty= 2m$, $m\in {\bf Z}$. It has  one branch when $\cos\pi\theta_\infty<-1$,  namely  $\theta_\infty=2m+1+i\vartheta$, and  when   $\cos\pi\theta_\infty>1$, namely $\theta_\infty=2m+i\vartheta$, $\vartheta>0$.  We divide the curve into seven portions: 

\vskip 0.2 cm 

[i] The half-line for $\cos\pi\theta_\infty<-9$ and  $p_{0x}<-2$. Here   $\sigma=1+2i\nu$, $\nu>0$. 

\vskip 0.2 cm 
[ii] The point  $(\cos\pi\theta_\infty,p_{0x})=(-9,-2)$. Here $\sigma=1$.  

\vskip 0.2 cm 
[iii] The segment of  line connecting $(-9,-2)$ to $(1,2)$, where  $-9<\cos\pi\theta_\infty<1$ and  $-2<p_{0x}<2$. Here $0<\sigma<1$. 

\vskip 0.2 cm 
[iv] The point $(\cos\pi\theta_\infty,p_{0x})=(1,2)$.  Here $\sigma=0$.  

\vskip 0.2 cm 
[v]  the segment of line connecting $(-1,-2)$ and $(1,2)$,   where $-1<\cos\pi\theta_\infty<1$ and $-2<p_{0x}<2$. Here  $0<\sigma<1$.

\vskip 0.2 cm 
[vi] The point $(\cos\pi\theta_\infty,p_{0x})=(-1,-2)$. Here  $\sigma=1$.

\vskip 0.2 cm 
[vii] The half-line  for $\cos\pi\theta_\infty>-1$ and $p_{0x}<-2$. Here  $\sigma=1+2i\nu$, $\nu>0$

\vskip 0.2 cm 
\noindent
The case of solutions (\ref{pro2y}), namely the case $
\sigma=1+2i\nu$, $\nu>0$, corresponds to the portions [i] and [vii] above. Included in this, is the  case of the quantum cohomology of $CP^2$, which   is on the portion [vii], for $(\cos\pi\theta_\infty,p_{0x})= (1,-7)$.

 \ble
Let $\sigma=1+2i\nu$, $\nu>0$, $\nu=(2/\pi)\ln G$. Let $
p_{0x}=p_{x1}=p_{01}=2\cos(\pi\sigma)=-2\cosh(2\pi\nu)\leq -2 
$. 
 The branch {\rm [i]$\cup$[ii]} of the cubic curve (\ref{cubic1}) for $(\cos\pi\theta_\infty,p_{0x})$ between $(-\infty,-\infty)$ and $(-9,-2)$ has equations: 
\be
\label{costeta1}
 \cos\pi\theta_\infty= -{1\over 2}e^{-3\pi\nu}\Bigl(3e^{5\pi\nu}+4e^{3\pi\nu}+3e^{\pi\nu}+e^{6\pi\nu}+3e^{4\pi\nu}+3e^{2\pi\nu}+1
\Bigr)
\ee
\be
\label{costeta1G}
=-{(G^2+G+1)(G^2-G+1)(G^8+2G^6+2G^2+1)\over 2G^6}
\ee
and:
$$
d= {i\over 2} \ln
\left\{
{4(G^4+1)^2\over (G^2-1)^2}~{16^{2i\nu}\Gamma\left({3-\theta_\infty\over 2}-i\nu\right)^2\Gamma\left({1+\theta_\infty\over 2}-i\nu\right)^2
\over \nu^2(2\nu+i(1-\theta_\infty))\Gamma(-i\nu)^4}
\right\}
$$
The branch {\rm [vii]$\cup$[vi]}  for $(\cos\pi\theta_\infty,p_{0x})$ between $(-1,-2)$ and $(+\infty,-\infty)$ has equations: 
\be
\label{costeta}
 \cos\pi\theta_\infty= -{1\over 2}e^{-3\pi\nu}\Bigl(3e^{5\pi\nu}+4e^{3\pi\nu}+3e^{\pi\nu}-e^{6\pi\nu}-3e^{4\pi\nu}-3e^{2\pi\nu}-1
\Bigr)
\ee
\be
\label{costetaG}
={(G^4-G^2+1)(G^8-2G^6-2G^2+1)\over 2G^6}
\ee
and:
\be
\label{dQC-Chazy}
d= {i\over 2} \ln
\left\{
{4(G^4+1)^2\over (G^2+1)^2}~{16^{2i\nu}\Gamma\left({3-\theta_\infty\over 2}-i\nu\right)^2\Gamma\left({1+\theta_\infty\over 2}-i\nu\right)^2
\over \nu^2(2\nu+i(1-\theta_\infty))\Gamma(-i\nu)^4}
\right\}
\ee
\ele

\vskip 0.2 cm 
\noindent
{\it Proof:} It is a matter of computation. To find (\ref{costeta1}) and (\ref{costeta}) substitute $p_{0x}=-2\cosh(2\pi\nu)$ in (\ref{cubic1}) and solve for $\cos\pi\theta_\infty$.  Substitute $\nu={2\over \pi}\ln G$ into (\ref{costeta1}) and (\ref{costeta})  and find (\ref{costeta1G}) and (\ref{costetaG}). These last are then substituted in (\ref{generald}), together with $p_{1x}=p_{01}=  p_{0x}=-2\cosh(2\pi\nu)$. Simple algebra gives the expression of $d$. 

\rightline{\qed}

%%%%%%%%%%%%%%%%%%%%%%%%%%%%%%%%%%

\subsection{Proof of Theorem \ref{theo2}}
\label{proofpro4}

Let $\mu=-1$ and $p_{0x}=p_{01}=p_{x1}=-7$, so that the branch $y(x)$ in (\ref{pro2y}) is associated to the quantum cohomology of $CP^2$. Then, from (\ref{nugen}) we immediately have $
\nu = {2\over \pi} \ln{1+\sqrt{5}\over 2}$. 
Substitution into (\ref{dQC-Chazy}) gives:
$$
 d= {i\over 2} \ln
\left\{
{4({\bf G}^4+1)^2\over ({\bf G}^2+1)^2}~{16^{2i\nu}\Gamma\left({5\over2}-i\nu\right)^2\Gamma\left(-{1\over 2}-i\nu\right)^2
\over \nu^2(2\nu+3i)\Gamma(-i\nu)^4}
\right\}
$$
Then, standard manipulation of $\Gamma$ functions proves (\ref{3vetrini}). The second part of the theorem follows from theorem \ref{theo1}. Please, note that we need already the result $\Im d= {\pi \nu\over 2}$ in order to compute the formula of the zeros. This is proved below.  

\rightline{\qed} 

\subsection{Proof of the Corollary} 
We recall the following convergent expansion 
\be
\label{ALK}
\ln(\Gamma(1+z)) = -\gamma~z-{1\over 2}\ln\left({\sin(\pi z)\over \pi z}\right)-\sum_{n=1}^\infty {\zeta(2n+1)z^{2n+1}\over 2n+1}, ~~~~~|z|<1
\ee
It can be applied to the factors $\ln(\Gamma(1-2i\nu))$ and $\ln(\Gamma(1-i\nu))$ in (\ref{3vetrini}),  provided that $|\nu|<1/2$. We also need to fix the determination 
$$
\ln\left({1-2i\nu\over 1+2i\nu}\right)= -2i\arccos{1\over \sqrt{1+4\nu^2}}
$$
with  arccos  positive for positive $\nu$ and zero for $\nu=0$.  With this in mind, we expand (\ref{3vetrini}) and find
$$
d=
{\pi\over 2}-8\nu\ln(2)+2\arccos{1\over \sqrt{1+4\nu^2}}+4 \sum_{n=1}^\infty {(-1)^n(1-4^n)\zeta(2 n+1)\over 2n+1}  ~ \nu^{2 n+1}~
+$$
$$+i\ln({\bf G}^2-1)~+k\pi,
$$ 
To conclude, observe that ${\bf G}^2-1={\bf G}$, by definition of golden ratio. Thus $\ln({\bf G}^2-1)=\ln {\bf G}= \pi\nu/2$. 

\rightline{\qed}

%%%%%%%%%%%%%%%%%%%%%%%%%%%%%%%%%%%%%%%%%%%%%%%%%%%5

\section{Proof of Proposition \ref{limitd}}
\label{prooflimitd}

 We distinguish two cases: 

\vskip 0.2 cm 

 1) Generic case   (\ref{generald}). To compute the expansion at $\nu=0$, we rewrite:
$$
 \Gamma\left({3\over 2}-\mu-i\nu\right)\Gamma\left( {1\over 2}+\mu-i\nu \right)
= {\pi\left({1\over 2}-\mu-i\nu\right)\over \sin\left[\pi\left(\mu+{1\over 2}+i\nu\right)\right]}{\Gamma\left({1\over 2}+\mu-i\nu\right)
\over \Gamma\left(\mu+{1\over 2}+i\nu\right)}
$$
The factor   ${\Gamma(z-x)/\Gamma(z+x)}$, $z=\mu+1/2$, $x=i\nu\to 0$, has Taylor series:
$$
{\Gamma(z-x)\over \Gamma(z+x)}=1-2\psi(z)x+2\psi(z)^2~x^2-\left({1\over 3}{d^2\psi (z) \over dz^2}+{4\over 3}\psi(z)^3\right)x^3+...
$$
The other factors in (\ref{generald}) are easily expanded around $\nu=0$, and the statement of Proposition \ref{limitd} follows. 

\vskip 0.2 cm 
2) Special cases (\ref{generald1}), (\ref{generald2}), (\ref{generald3}), (\ref{generald4}). We consider the case of  (\ref{generald1}), the others being analogous.  Manipulation of $\Gamma$ functions gives: 
$$
{\nu^2\Gamma(m)^2\Gamma(2i\nu+1-m)^2(p_{x1}-2)
\over
16^{2i\nu}\Gamma(1+i\nu)^4
}
=
{2-p_{x1}\over 4} \left[
1+\sum_{k=1}^{m-1} {4i\nu\over k}+O(\nu^2)
\right] e^{-8i\nu\ln 2}{\Gamma(1+2i\nu)^2\over \Gamma(1+i\nu)^4}
$$
Recall that $\sum_{k=1}^{m-1} {1\over k}=\psi(m)+\gamma$. Also, expand the exponential $ e^{-8i\nu\ln 2}$. We find: 
$$
d= -{i\over 2}\ln\left({2-p_{x1}\over 4}\right)-{i\over 2}\ln\left(1+4i\nu(\psi(m)+\gamma-2\ln2)+O(\nu^2)\right)-i\ln{\Gamma(1+2i\nu)\over \Gamma(1+i\nu)^2}
$$
The logarithm of $\Gamma$ functions is $O(\nu^2)$, as it is clear from the expansion (\ref{ALK}). The statement of Proposition \ref{limitd} follows. 

\rightline{\qed}

\vskip 0.2 cm 

Proposition \ref{limitd} fails in the singularities   $2\mu=2m+1$, $m\in{\bf Z}$, of the expansion of $d$ for $\nu \to 0$. 

%%%%%%%%%%%%%%%%%%%%%%%%%%%%%%%%%%%%%%%%%%%%%%%%%%%%%%%%%%

\section{Proof of Proposition \ref{chazyd}}
\label{proofchazyd}

\ble
\label{tetanu} Let ${\bf G}$ be the golden ratio. 
On the segment of the curve (\ref{cubic1}) between $(-1,-2)$ and $(1,-7)$, $\theta_\infty$ is real and $
0\leq \nu\leq {2\over \pi} \ln {\bf G}$. 
As a function of $\nu$, $\theta_\infty$ is: 
\be
\label{teta}
\theta_\infty = \pm{1\over \pi}\hbox{\rm arg}\left(c_1(\nu)-{i\over 2}\sqrt{c_2(\nu)}\right)+2m+1,~~~~~m\in{\bf Z}
\ee
where $c_1(\nu)$ and $c_2(\nu)$ are given by the expressions (\ref{starrr}). In particular, $c_2(\nu) \geq 0$ and the square root is the positive one. Moreover  $
\left|c_1(\nu)-{i\over 2}\sqrt{c_2(\nu)}\right|=1
$. 
The argument is chosen such that $
-\pi\leq \hbox{\rm arg} \left(c_1(\nu)-{i\over 2}\sqrt{c_2(\nu)}\right)\leq 0
$.

\ele

\noindent{\it Proof:} It is clear that $\theta_\infty $ is real, because $-1\leq \cos\pi\theta_\infty\leq 1$.  Solving (\ref{costeta})  for $\theta_\infty$ we find: 
$$
\theta_\infty =  -{i\over \pi}\ln\left(c_1(\nu)\pm{i\over 2}\sqrt{c_2(\nu)}\right)+2m+1,~~~~~m\in{\bf Z}
$$
The discriminant $c_2(\nu)$ is positive for $-{2\over \pi} \ln {\bf G}<\nu<{2\over \pi} \ln {\bf G}$ and simple computation shows that in this case $|c_1(\nu)-{i\over 2}\sqrt{c_2(\nu)}|=1$. Therefore:
$$
\theta_\infty =
 {1\over \pi}\hbox{\rm arg}\left(c_1(\nu)\pm{i\over 2}\sqrt{c_2(\nu)}\right)+2m+1=
\pm{1\over \pi}\hbox{\rm arg}\left(c_1(\nu)-{i\over 2}\sqrt{c_2(\nu)}\right)+2m+1
$$

\rightline{\qed}

\vskip 0.2 cm 
Due to the symmetries of $PVI_\mu$, one can restrict to the case $-1\leq \Re \mu<0$, i.e. $-2\leq \Re(\theta_\infty)<0$, and take  
\be
\label{raffreddore}
 \theta_\infty = -1+{1\over \pi}\hbox{\rm arg}\left(c_1(\nu)-{i\over 2}\sqrt{c_2(\nu)}\right)
\ee
In particular, for  $\nu={2\over \pi} \ln {\bf G}$ ($p_{0x}=p_{x1}=p_{01}=-7$):
$$ 
 c_1(\nu)-{i\over 2}\sqrt{c_2(\nu)}=-1~~~\Longrightarrow~~~\theta_\infty= -2 ~~~\hbox{(Quantum Cohomology of $CP^2$, $\mu=-1$)}
$$
For $\nu =0$ ($p_{0x}=p_{x1}=p_{01}=-2$):
$$ 
 c_1(\nu)-{i\over 2}\sqrt{c_2(\nu)}=1~~~\Longrightarrow~~~\theta_\infty= -1,~~~ \hbox{(Chazy solutions of \cite{MazzChazy}, $\mu=-{1\over 2}$)}
$$

\vskip 0.5 cm
\noindent
{\it Proof of Proposition \ref{chazyd}:} Observe that (\ref{raffreddore}) implies:
$$
  \theta_\infty= -1+2\omega(\nu),~~~~~~~\nu\to 0
$$
$$
2\omega(\nu):=-{\pi\sqrt{3}}~\nu^2-{\sqrt{3}\pi^3\over 4}~\nu^4-{17\over 120}\pi^5\sqrt{3}~\nu^6+O(\nu^8)
$$
 Substitute this into (\ref{dQC-Chazy}) and use the identities $\Gamma(z+1)=z\Gamma(z)$, $\Gamma(1-z)=\pi\Gamma^{-1}(z)\sin^{-1}(\pi z)$ to simplify. We obtain: 
$$
 d={i\over 2} 
\ln
\left(
-{
(G^4+1)^2 16^{2i\nu}\nu^2\Gamma^2(1-\omega-i\nu)\Gamma^2(1+\omega-i\nu)^2
\over 
(G^2+1)^2(i\nu-\omega)^2\Gamma^4(1-i\nu)
}
\right),~~~~~G=e^{\pi\nu\over 2}.
$$
When $\nu\to 0$, then $|i\nu\pm \omega|<1$. Thus, one can use the expansion (\ref{ALK}), applied to $\Gamma(1\pm\omega-i\nu)$ and $\Gamma(1-i\nu)$. 
This gives the expansion of  Proposition \ref{chazyd}.

\rightline{\qed}

%%%%%%%%%%%%%%%%%%%%%%%%%%%%%%%%%%%%%%%%%%%%%%%%%%%%%%%%%%%%%%
\section{Appendix: Coefficients $A_{nm}$}\label{coeff}

\noindent
Order 3:
$$
A_{33}={3\over 2^{16}}{(2\mu-1-2i\nu)^6 \over \nu^6}
$$
$$
A_{32}={i\over 2^{15}}{(-2\mu+2i\nu+1)^4
\over 
(-\nu+i)^2\nu^6
}\times
$$
$$
\times
(-128\nu^5+288i\nu^4+192\nu^3+8i(2\mu+1)(2\mu-3)\nu^2-9(i-2\nu)(2\mu-1)^2
)
$$
$$
A_{31}={(-2\mu+2i\nu+1)^2\over 2^{16} \nu^6(i-\nu)}(2048\nu^7-560i\nu^4-1616\nu^5-1440i\mu^3-1376i\nu^2\mu-360i\mu-840\nu^3-3360\nu^3\mu^2+
$$
$$+3360\nu^3\mu+128i\nu^2\mu^4+720i\mu^4-256i\nu^2\mu^3-2048i\nu^4\mu^2+1504i\nu^2\mu^2-45
\nu-720\nu\mu^4+1440\nu\mu^3+
$$
$$
-1080\nu\mu^2+360\nu\mu+336i\nu^2+2048i\nu^4\mu+45i-3712i\nu^6+1080i\mu^2)
$$
$$
A_{30}=-{(2\mu-1)^2+4\nu^2\over 2^{14} (\nu^2+1)^2\nu^6}
(320\nu^4\mu^4+528\nu^2\mu^4+240\mu^4-640\nu^4\mu^3-1056\nu^2\mu^3-480\mu^3+$$
$$
+
448\nu^6\mu^2+1328\nu^4\mu^2+1224\nu^2\mu^2+360\mu^2-448\nu^6\mu-1008\nu^4\mu-696\nu^2\mu +
$$
$$
-120\mu-128\nu^8-144\nu^6+104\nu^4+141\nu^2+15)
$$

\noindent
Order 4: 
$$
A_{44}= -{ (-2\mu+2i\nu+1)^8 \over 2^{20} \nu^8}
$$
$$
A_{43}=-{i(-2\mu+2i\nu+1)^6\over 2^{19}(\nu-i)^2\nu^8}\times
$$
$$
\times(-72\nu^5+168i\nu^4+120\nu^3+12i\nu^2\mu^2-12i\nu^2\mu-21i\nu^2+8\nu+32\nu\mu^2-32\nu\mu-4i-16i\mu^2+16i\mu)
$$
$$
A_{42}={ (-2\mu+2i\nu+1)^4 \over 2^{18}(\nu-i)^2\nu^8}(7+1120i\nu^5-2976i\nu^7-224\mu^3+306i\nu^3-56\mu+336i\nu\mu^2-112i\nu\mu+
$$
$$-448i\mu^3\nu+1264i\nu^3\mu^2+512i\nu^5\mu+32i\mu^4\nu^3-64i\nu^3\mu^3-512i\mu^2\nu^5+224i\mu^4\nu+
$$
$$-1232i\nu^3\mu-288\nu^2\mu-64\nu^2\mu^4
+
128\nu^2\mu^3-1440\nu^4\mu^2+224\nu^2\mu^2+
$$
$$
+1440\nu^4\mu+14i\nu+112\mu^4-3040\nu^6+1024\nu^8+168\mu^2+76\nu^2-200\nu^4)
$$
$$
A_{41}= { i  (-2\mu+2i\nu+1)^2\over 3~2^{19} (i-\nu)(\nu+i)^2\nu^8 }(-84-36560\nu^7-6504i\nu^5-53248i\nu^9+26880\nu^8\mu+40704\nu^6\mu^3+
$$
$$-20352\mu^4\nu^6-8256\nu^2\mu^6-4608\nu^4\mu^5-26880\mu^2\nu^8+24768\nu^2\mu^5+1536\mu^6\nu^4+13440\mu^3+5376i\mu^6\nu+
$$
$$
+51264i\nu^7\mu+3072i\nu^7\mu^4+42768i\mu^4\nu^3-48192i\nu^7\mu^2+1008\mu+59904\nu^6\mu+16452i\mu^2\nu^3-13440\mu^3\nu
+$$
$$-34656i\mu^3\nu^3+23616i\nu^5\mu+24576i\nu^9\mu+10176i\nu^3\mu^6-19200i\nu^5\mu^3-16128i\mu^5\nu-80256\nu^6\mu^2+
$$
$$-37632\nu^4\mu^4+82944\nu^4\mu^3+12492\nu^2\mu-52848\nu^2\mu^4+64416\nu^2\mu^3-86688\nu^4\mu^2-40572\nu^2
\mu^2+44448\nu^4\mu+
$$
$$
+447i\nu^3+84i\nu-23552i\nu^{11}-20160\mu^4-12552\nu^6+5040i\mu^2\nu-24576i\nu^9\mu^2-4212i\mu\nu^3+20160i\mu^4\nu+
$$
$$-1008i\mu\nu+9600i\nu^5\mu^4-6144i\nu^7\mu^3-30528i\nu^3\mu^5+11968\nu^8-14016i\nu^5
\mu^2+16384\nu^12+
$$
$$
-5376\mu^6+16128\mu^5+33920\nu^10-5040\mu^2-1497\nu^2-8448\nu^4
)
$$
$$
A_{40}=-{  (2\mu^2-1)^2-4\nu^2 \over 2^{19} (\nu^2+1)^2\nu^8}(
35-10048\nu^8\mu-25216\nu^6\mu^3+12608\mu^4\nu^6+5760\nu^2\mu^6-12096\nu^4\mu^5+
$$
$$
+10048\mu^2\nu^8-17280\nu^2\mu^5+4032\mu^6\nu^4-5600\mu^3-420\mu-24864\nu^6\mu+37472\nu^6\mu^2+35472\nu^4\mu^4+
$$
$$
-50784\nu^4\mu^3-5720\nu^2\mu+30880\nu^2\mu^4-32960\nu^2\mu^3+43844\nu^4\mu^2+19320\nu^2\mu^2-20468\nu^4\mu+8400\mu^4+
$$
$$
+2804\nu^6-2736\nu^8+2240\mu^6-6720\mu^5-2624\nu^{10}+2100\mu^2+670\nu^2+3719\nu^4
)
$$

%%%%%%%%%%%%%%%%%%%%%%%%%%

\vskip 0.5 cm
\noindent 
{\bf Acknowledgments:}   I thank Korea Institute for Advanced Study for providing computing resources
(Abacus System) for this work. This work was initiated as a result of fruitfull discussions with M. Mazzocco.  

%%%%%%%%%%%%%%%%%%%%%%%%%%%%%%%%%%%%%%%%%%%%%%%%%%%%%%%%%%%%%%%%%%%%%%%%%%%%%%%%%%%%%%%%%%%%%%%%%

\end{document}